%% file: br
\documentclass[12pt]{amsart}
\usepackage{amssymb, amsmath, epsfig,graphics}

\def\mapdownright#1{\Big\downarrow\rlap{$\vcenter{\hbox{$\scriptstyle#1$}}$}}

\def\mapse#1{\rlap{$\vcenter{\hbox{$\scriptstyle#1$}}$}\,\,\searrow}

\newtheorem{theore}{Theorem}[section]
\newtheorem{propositio}[theore]{Proposition}
\newtheorem{lemm}[theore]{Lemma}
\newtheorem{remar}[theore]{Remark}
\renewcommand{\proof}{\vspace{.05in}
                    \noindent {\sc Proof} \hspace{.05in}}
\newcommand{\ethrm}{\hspace*{\fill}
                      $\Box$
                      \vspace{.1in}}

\newcommand{\aut}{{\mathop{\rm Aut}\nolimits}}
\newcommand{\ra}{\rightarrow}

\newcommand{\pic}{\mathop{\rm Pic\,}\nolimits}

\renewcommand{\phi}{\varphi}
\newcommand{\SL}{{\rm SL}}
\newcommand{\GL}{{\rm GL}}

\renewcommand{\Box}{\square}
\renewcommand{\tilde}{\widetilde}

\newcommand{\al}{\alpha}

\newcommand{\la}{\lambda}

\newcommand{\M}{\mathcal M}

\newcommand{\A}{\mathbb A}
\newcommand{\e}{_{\text{\'et}}}

\newcommand{\h}{{\mathfrak h}}

\newcommand{\X}{\mathcal X}

\newcommand{\Y}{\mathcal Y}
\newcommand{\E}{\mathcal E}
\newcommand{\C}{\mathbb C}
\newcommand{\Z}{\mathbb Z}
\renewcommand{\O}{\mathcal O}
\newcommand{\U}{\mathcal U}

\renewcommand{\P}{\mathbb P}

\newcommand{\R}{\mathbb R}
\newcommand{\sH}{\mathcal H}

\newcommand{\color}[6]{}

\newcommand{\db}[1]{D^b(#1)}
\newcommand{\daut}[1]{{\rm Auteq}(#1)}
\newcommand{\rhom}{{\bf R}{\mathcal H}\!\text{\em om}}
\newcommand{\tensor}{\otimes^{{\bf L}}}
\newcommand{\push}{{\bf R}}
\newcommand{\pull}{{\bf L}}
\newcommand{\cone}{\mathop{\rm Cone}\nolimits}

\begin{document}
\begin{center}
{\LARGE Artin group actions on derived categories of threefolds} 

\vspace{0.2in}

{\large Bal\'azs Szendr\H oi}

\vspace{0.15in}

{\large October 2002}

\end{center}

\vspace{0.15in}

{\small\begin{center} {\sc abstract} \end{center}
{\leftskip=30pt \rightskip=30pt
Motivated by the enhanced gauge symmetry phenomenon of the physics literature and 
mirror symmetry, this paper constructs an action of an Artin group on the derived 
category of coherent sheaves of a smooth quasiprojective threefold containing 
a configuration of ruled surfaces described by a finite type Dynkin diagram. 
The action extends over deformations of the threefold via a compatible 
action of the corresponding reflection group on the base of its 
deformation space. All finite type Dynkin diagrams are realized.  
\par}
\vspace{0.15in}
\begin{center}
\end{center}} 
\vspace{0.15in}

\section*{Introduction} 

The purpose of this work is to construct actions of finite type Artin groups 
on derived categories of coherent sheaves of complex threefolds. 
The construction is motivated by a correspondence between Calabi--Yau threefolds 
containing ruled surfaces and Lie algebras, which arises in 
Type II string theory; see~\cite{kmp},~\cite{akm},~\cite{mg} and references in 
these works, as well as~\cite{c..v} which also considers some of the geometries
studied in this paper. 
I~explain the connection further in~\cite{eg}; suffice it to say here
that the threefolds I~consider are of the most simple kind for which the 
physics correspondence works, and the main theorem of this paper 
says that in these cases the 
derived category of coherent sheaves of the threefold is acted on by an Artin 
group which covers the Weyl group of the corresponding Lie algebra. 

The main result of this paper can also be viewed in the context of homological 
mirror symmetry, representing a generalization of a result of Seidel and 
Thomas~\cite{st}. 
They construct representations of the classical (Type~$A$) braid group on 
derived categories of coherent sheaves of a much larger class of varieties 
than those considered here. However, their construction is more algebraic in 
flavour, whereas the Artin group actions in this paper are governed in a very
precise geometric way by deformation theory. As explained in~\cite{eg}, 
the two constructions coincide in dimension two; the braid group actions on derived 
categories of threefolds obtained in this paper are new even in the Type~$A$ case.

Autoequivalences of derived categories for threefolds containing ruled 
surfaces were first constructed by Horja in~\cite{horja1}--\cite{horja2}.  
In~\cite{sz} it was observed that these equivalences are essentially 
given by classical correspondences (structure sheaves of 
subschemes in the product), and also that they deform to  
derived equivalences given by flops first found by Bondal and Orlov~\cite{bo}.
The proof of the braid relations uses in an essential way both of these facts. 
 
In certain cases, the Artin group acts faithfully on the derived category. 
The proof of this statement will be reduced to the 
injectivity statement of Seidel--Thomas~\cite{st} for Type~$A$, using a 
hyperplane section argument. 

\vspace{0.1in}

\noindent{\bf Structure of the paper} Section~\ref{sec!dyn} deals with 
Dynkin diagrams, reflection groups and braid groups. 
In Section~\ref{sec!23} I 
first recall some results about resolutions of Kleinian surface singularities, 
and then turn to the construction of certain quasiprojective 
threefolds and their deformations. 
Section~\ref{sec!dq} discusses generalities about families of Fourier--Mukai 
functors. Section~\ref{sec!main} contains the main results. 
Families of Fourier--Mukai functors constructed in Section~\ref{sub!rel} 
are shown to satisfy braid relations in Section~\ref{sub!main}, where faithfulness
is also proved in some cases. 
The paper is concluded in Section~\ref{sub!proj} 
by a brief discussion of the projective case. 

\vspace{0.1in}

\noindent{\bf Conventions} A {\em smooth family} means a smooth morphism 
$e\colon \X\ra S$ of smooth varieties over~$\C$ with~$\X$ quasiprojective over $S$. 
The base $S$ will always be very simple in this paper, typically affine space
$\A^r$ or an open set thereof. For a brief period in Section~\ref{sec!23}, $r$ can 
be infinite, but this will cause no complications.  
The fiber~$e^{-1}(s)$ over~$s\in S$ will be denoted by~$X_s$.
By definition a {\it Dynkin diagram}~$\Delta$ means an irreducible diagram  
of finite type $A_n\ldots G_2$. Nodes of~$\Delta$ will be denoted $i,j,\ldots$;
for~$i\neq j$, $m_{ij}=m_{ji}\in\{2,3,4,6\}$ is the label associated to the 
pair of nodes~$(i,j)$. As usual, the pair~$(i,j)$ is said to span an edge 
if $m_{ij}>2$. The diagram~$\Delta$ is simply laced (type~$ADE$)
if~$m_{ij}\in\{2,3\}$.

\vspace{0.1in}

\noindent{\bf Acknowledgement} I~thank Roger Carter, Mark Gross, Sheldon Katz and  
J\'anos Koll\'ar for helpful remarks and correspondence. I~especially thank
Ian Grojnowski for his help in clarifying old ideas and coming up with new ones, 
during many conversations on a joint project which grew out of the present paper. 
Finally I~wish to acknowledge that I~learned the importance of structure sheaves 
from Tom Bridgeland's beautiful paper~\cite{bridgeland_flops}. 

\section{Dynkin diagrams and Artin groups}\label{sec!dyn}

\subsection{The reflection group and the Artin group}\label{braid}

Take an arbitrary Dynkin diagram~$\Delta$ with~$n$ nodes.  
Let $\Sigma_\Delta\subset\h_{\Delta,\R}$ be the corresponding root system, 
where $(\h_{\Delta,\R},\langle,\rangle)$ is a Euclidean inner product space of 
dimension~$n$. Fix sets of simple and positive roots 
\[\Sigma_\Delta^0=\{\la_1, \ldots, \la_n\}\subset\Sigma_\Delta^+\subset\Sigma_\Delta.\] 
The reflections $r_i\colon \h_{\Delta,\R}\ra\h_{\Delta,\R}$ $(1\leq i\leq n)$
defined by the simple roots generate a finite reflection group 
\[W_\Delta=\langle r_i\rangle<\GL(\h_{\Delta,\R}).\] As an abstract group, 
\[ W_\Delta\cong\Big\langle r_i\colon i\in{\rm Nodes}(\Delta)\Big\rangle\Big/\Big\langle r_i^2 = 1, (r_i r_j)^{m_{ij}}= 1 \Big\rangle\]
with one relation for every node~$i$ and one for every pair of different nodes~$(i,j)$
with label~$m_{ij}$. The set of reflections in~$W_\Delta$ is in one-to-one 
correspondence with the set $\Sigma_\Delta^+$. The group~$W_\Delta$ also acts on 
the complex vector space $\h_\Delta=\h_{\Delta,\R}\otimes\C$; for a reflection 
$w\in W_\Delta$, let $\Pi_w\subset\h_\Delta$ denote the fixed hyperplane of~$w$. 

Define the {\it Artin group} (also called generalized braid group)~$B_\Delta$ 
by generators and relations as 
\begin{equation}\label{braidrelations} B_\Delta  = \Big\langle R_i\colon i\in{\rm Nodes}(\Delta)\Big\rangle\Big/\Big\langle \underbrace{R_i R_j\ldots}_{m_{ij}} = \underbrace{R_j R_i\ldots}_{m_{ij}}\Big\rangle 
\end{equation}
with one relation for every pair of different nodes~$(i,j)$ of~$\Delta$, the 
{\em braid relation}. 
There is a group homomorphism $B_\Delta\rightarrow W_\Delta$ sending 
$R_i$ to~$r_i$. 

\subsection{Quotiens of Dynkin diagrams}\label{subsec!qdyn}

Let~$\Delta$ be a simply laced Dynkin diagram and~$A$ a non-trivial subgroup of 
its automorphism group~$\aut(\Delta)$, excluding the case 
$(\Delta, A)=(A_{2n}, \Z/2)$. 
Then~$A$ permutes the set of simple roots
in~$\h_\Delta$ which forms a basis; hence~$A$ also acts on~$\h_\Delta$.
Let~$\h_\Xi=(\h_\Delta)^A$ and let~$\Sigma_\Xi=\Sigma_\Delta\cap \h_\Xi$ 
be the set of invariant roots. It is well known that, as the notation suggests,  
$\Sigma_\Xi$ is a root system for a Dynkin diagram~$\Xi$, the 
``quotient'' diagram~$\Delta/ A$. For~$\Delta=A_{2n-1}, D_n, E_6$,~$\Xi$ is 
the non-simply laced diagram of type~$C_n, B_{n-1}, F_4$ respectively; 
in the case~$\Delta=D_4$,~$\Xi$ is either~$G_2$  
or~$C_3$ according to whether~$A$ acts transitively on the outer nodes
of~$D_4$ or not. The set of nodes of~$\Xi$ is in one-to-one correspondence with 
the set of orbits of nodes of~$\Delta$ under the action of~$A$; 
there is a corresponding set of simple and positive roots 
\[\Sigma_\Xi^0=\{\mu_1, \ldots \}\subset\Sigma_\Xi^+\subset\Sigma_\Xi\] and 
a reflection group~\[W_\Xi=\langle\rho_i\colon i\in{\rm Nodes}(\Xi)\rangle\]
acting on~$\h_\Xi$. Finally~$\Xi$ also defines an Artin group~$B_\Xi$. 

\renewcommand{\labelenumi}{(\roman{enumi})}
\begin{lemm} The group~$A$ acts naturally on the Artin group~$B_\Delta$ and the 
reflection group~$W_\Delta$ equivariantly with respect to the map 
$B_\Delta\rightarrow W_\Delta$.
The fixed subgroups are isomorphic to~$B_\Xi$ and~$W_\Xi$ respectively.
\label{maps_braidgp}\end{lemm}
\proof For~$a\in A$, the action is defined on generators of~$B_\Delta$ by
$R_i\mapsto R_{a(i)}$. This action clearly leaves the relations invariant and
descends to an action on~$W_\Delta$. If~$\{k_j\}$ is an~$A$-orbit of nodes
of~$\Delta$ corresponding to a node~$k$ of~$\Xi$, 
then~$R_{k_j}\in B_\Delta$ commute by the braid relations and their product 
$R_k=\prod R_{k_j}$ is invariant under the action. It is an easy check to show 
that the~$R_k$ satisfy the braid relations of the group~$B_\Xi$. 
By~\cite[Corollary 4.4]{michel}, these elements generate the fixed subgroup 
and they do not satisfy any further relations. 
\ethrm 

\begin{remar} \rm\label{acts} The proof also shows that 
the action of~$A$ on the reflection group~$W_\Delta$ is simply the 
conjugation action of~$A<\GL(\h_\Delta)$ on~$W_\Delta<\GL(\h_\Delta)$. 
\end{remar}

\section{Surfaces, threefolds and deformations}\label{sec!23}

\subsection{Finite subgroups of~$\SL(2,\C)$}\label{mckay}

Fix a finite subgroup~$\Gamma$ of~$\SL(2,\C)$ together with its canonical 
two-dimensional representation~$\rho_c$. Following McKay, 
consider the diagram~$\tilde\Delta$ consisting of a node for every 
irrep (irreducible representation)~$\rho_j$ of~$\Gamma$ and an edge 
between irreps~$\rho_j$,~$\rho_k$ whenever~$\rho_j$ is a direct 
summand of~$\rho_k\otimes\rho_c$. It is well known that this defines a 
symmetric relation and 
$\tilde\Delta$ is an affine Dynkin diagram 
of type~$\tilde A_n, \tilde D_n$ or~$\tilde E_n$ with 
distinguished affine node corresponding 
to the trivial rep~$\rho_0$. Let~$\Delta=\tilde\Delta\setminus\{\rho_0\}$ 
be the corresponding finite diagram. 

\begin{lemm} There exists an exact sequence of groups
\begin{equation}\label{seq_of_Gamma} 1 \ra C_\Gamma \stackrel{\gamma}{\ra} N_\Gamma \stackrel{\delta}{\ra} \aut(\Delta) \ra 1, \end{equation}
where~$N_\Gamma=N_{\GL(2,\C)}(\Gamma)/\Gamma$ and 
$C_\Gamma=C_{\GL(2,\C)}(\Gamma)/Z_\Gamma$ are the normalizer of~$\Gamma$ in 
$\GL(2,\C)$ modulo~$\Gamma$ and the centralizer
of~$\Gamma$ modulo the center~$Z_\Gamma$ of~$\Gamma$ respectively,  
and~$\aut(\Delta)$ is the automorphism group of the diagram~$\Delta$. 
\label{groups_lemma}\end{lemm}
\proof 
For an irrep~$\rho\colon\Gamma\ra \GL(V)$ and 
an element~$g\in N_{\GL(2,\C)}(\Gamma)$, define a new irrep~$\rho^g$ of 
$\Gamma$ by~$\rho^g(h)=\rho(g^{-1}h g)$. The isomorphism
class of the irrep~$\rho^g$ only 
depends on the class of~$g$ in~$N_\Gamma$. For all 
$g\in N_{\GL(2,\C)}(\Gamma)$,~$\rho_0^g$ is isomorphic to~$\rho_0$ and 
$\rho_c^g$ to~$\rho_c$, so the diagram~$\Delta$ is mapped to itself by the action
of~$g$. So~$\rho\mapsto \rho^g$ defines a map 
$\delta\colon N_\Gamma \ra \aut(\Delta)$. The proof of the surjectivity of this map 
as well as the computation of its kernel are easy on a case-by-case basis. 
\ethrm

\subsection{The surface~$Y$}
For a finite subgroup~$\Gamma<\SL(2,\C)$, let~$g\colon Y\ra \C^2/\Gamma$ be the 
minimal resolution of the Kleinian quotient 
singularity with exceptional locus~$E=\cup_{j=1}^n E_j$. 
The incidence graph of the components of~$E$ can be identified 
with the (simply laced) McKay diagram~$\Delta$ defined by~$\Gamma$ in~\ref{mckay}; 
fix such an identification. 

\renewcommand{\labelenumi}{(\roman{enumi})}
\begin{propositio} There is an injection
$j\colon N_\Gamma \hookrightarrow \aut(Y)$. 
The composite \[N_\Gamma\ra \aut(Y)\ra\aut(\Delta),\] where the latter is the map 
given by permuting exceptional divisors, coincides with the map~$\delta$ of 
Lemma~\ref{groups_lemma}. In particular, an element of~$N_\Gamma$ 
fixes all the exceptional divisors if and only if it is in~$C_\Gamma$.
\label{auts_of_ALE}\end{propositio}
\proof By definition, an element~$h\in N_{\GL(2,\C)}(\Gamma)$ 
induces an automorphism of~$\C^2$ normalizing the action of~$\Gamma$; 
hence this automorphism descends to 
the quotient~$\C^2/\Gamma$ and only depends on the class of~$h$ in 
$N_\Gamma$. The resolution~$f$ is the unique minimal model of~$\C^2/\Gamma$, hence 
every automorphism of~$\C^2/\Gamma$ lifts to a 
unique automorphism of~$Y$. This defines the injection 
$j\colon N_\Gamma \ra \aut(Y)$. The last statement follows from an explicit 
computation on the resolution. 
\ethrm 

Next I~collect information about the cohomology and deformations 
of the surface~$Y$. The first statement is well known. 

\begin{propositio}\label{coh_of_Y}
The second cohomology~$H^2(Y, \Z)$ of the surface~$Y$ is a free~$\Z$-module
of rank~$n$, with dual~$H^2_c(Y,\Z)$ which has a natural~$\Z$-basis consisting of 
the classes~$\{[E_j]\}$ of the exceptional divisors. 
Every exceptional divisor~$E_j\subset Y$ gives rise to a reflection
\[\omega\mapsto\omega + \left([E_j]\cdot\omega\right)c_1(\O_Y(E_j)).\]
on~$H^2(Y,\C)$. These reflections generate a finite reflection group. 
\end{propositio}\ethrm

\begin{remar} \rm The data in this proposition can be identified with the 
Dynkin diagram data as follows. Mapping the class~$[E_j]\in H^2_c(Y,\C)$ to the 
simple root~$\la_j$ gives an isomorphism between the lattice~$H_c^2(Y, \Z)$ 
and the root lattice~$\Z\langle\la_1, \ldots, \la_r\rangle$. Dually, 
$H^2(Y, \C)\cong\h_\Delta$ and the reflection defined by~$E_j$ is the 
reflection~$r_j$ associated to the simple root~$\la_j$. Hence the reflection group 
in (ii) is isomorphic to the reflection group~$W_\Delta$. Note also that the group 
$\aut(\Delta)$ acts both on $H^2(Y,\C)$ (by acting on a basis) and on $\h_\Delta$
(as defined in Section~\ref{subsec!qdyn}) and these actions are obviously 
compatible.  
\end{remar}

Recall the hyperplanes~$\Pi_w\subset\h_\Delta$ which are the fixed loci 
of reflections of~$W_\Delta$; in particular, these include the fixed hyperplanes 
$\Pi_{r_j}= \{\omega\in\h_\Delta\, | \langle\omega,\la_i\rangle=0\}$ of~$r_j$. 

\begin{propositio} \ \label{defs_of_Y}\begin{enumerate}
\item The universal deformation space of~$Y$ is a smooth family 
$d\colon\Y\ra Z$ with central fiber~$d^{-1}(0)\cong Y$. There is a simultaneous 
contraction~$G\colon \Y\ra\bar \Y$ over~$Z$ with 
central fiber~$g\colon Y\ra \bar Y$. More generally, for any subset 
$I=\{E_{j_i}\}$ of the exceptional curves, there exists
a contraction~$G_I\colon \Y\ra\bar\Y_I$ over~$Z$ 
contracting curves in~$I$ in the central fiber and giving an isomorphism on 
fibers~$Y_t$ to which none of the curves in~$I$ deform. 
\item Choosing a generator of the relative canonical bundle 
$\omega_{\Y/Z}$ over~$Z$ gives rise to a period map
\[ \phi\colon Z\longrightarrow H^2(Y, \C)\cong\h_\Delta\]
which is an isomorphism. Different choices of the generator
give rise to the same isomorphism up to multiplication by a nonzero constant 
on~$\h_\Delta$. 
\item The group~$N_\Gamma$ acts naturally on the base~$Z$ and compatibly on the 
total space of the family~$\Y\ra Z$ by automorphisms. There is an induced 
action of~$N_\Gamma$ on~$\h_\Delta$ via $\phi$. This action factors via the 
morphism 
\[(\delta, \det)\colon N_\Gamma\ra \aut(\Delta)\times \C^*\]
with~$\C^*$ acting on~$\h_\Delta$ with weight one. 
\item The action of the reflection group~$W_\Delta$ on~$\h_\Delta$ induces, via the 
period map~$\phi$, a~$W_\Delta$-action on~$Z$. 
This extends compatibly to an action on the total space  
$\bar\Y\ra Z$ by automorphisms. More generally, if~$w\in W$ is a reflection and 
$I$ denotes the set of nodes corresponding to positive roots mapped to negative
roots by~$w$, then there is a diagram 
\[ \begin{array}{ccc}\bar\Y_I& \stackrel{\sim}{\longrightarrow}&\bar\Y_I\\
\Big\downarrow &&\Big\downarrow \\
Z& \stackrel{w}{\longrightarrow}&Z.\end{array}\]
There is no~$W_\Delta$-action by automorphisms on the total space~$\Y\ra Z$, 
but for~$w\in W_\Delta$ and~$t\in Z$ the fibers~$g^{-1}(t)$,~$g^{-1}(tw)$ are 
isomorphic. 
\item For a point~$s\in \h_\Delta$, the fiber~$g^{-1}\phi^{-1}(s)=Y_s$ contains 
projective curves if and only if~$s\in\cup_{w\in\Sigma_\Delta}\Pi_w$; 
otherwise it is affine. For~$w\in\Sigma_\Delta^+$ corresponding to the positive
root~$\la$, write~$\la=\sum_j a_j\la_j$. Then 
$s\in \Pi_w$ if and only if~$Y_s$ contains a rational (possibly reducible) curve
which is a flat deformation of the effective rational curve~$\sum_j a_jE_j$ on~$Y$.
If~$s\in\Pi_w\setminus\cup_{w'\neq w}\Pi_{w'}$ then this curve is smooth and 
it is the unique projective curve in~$Y_s$. 
\end{enumerate} 
\end{propositio}
\proof (i) is well known. (ii) is proved for example in~\cite[Section 4]{nisb}. 
For (iii), note that by universality, $\aut(Y)$ and hence its subgroup $N_\Gamma$ 
act on the family~$\Y\ra Z$, and in particular on its base $Z$. 
By~\cite[Proposition 8.6ii]{slo}, the determinant acts by weight one in the induced
action on $\h_\Delta$.
The action of the reflection group~$W_\Delta$ on $Z$ and its properties 
stated in~(iv) go back to Brieskorn, see e.g.~\cite{slo}. 
Finally (v) is spelled out in~\cite[Theorem 1]{km}. 
\ethrm 
\begin{remar} \rm Using (ii), I~will identify~$Z$ with~$\h_\Delta$ everywhere 
below. The space~$\h_\Delta$ carries actions of both~$N_\Gamma$ and~$W_\Delta$.  
\end{remar}

\subsection{The group cohomology class~$\al$ and related constructions}\label{aux} 

Let~$B$ be a smooth, not necessarily projective curve over~$\C$. 
Fix a cohomology class~\[\al\in H^1(B\e, N_\Gamma)\] in the 
nonabelian group cohomology set~$H^1(B\e, N_\Gamma)$ in the \'etale topology
(compare~\cite[p.122]{mi}). 
All constructions below depend on this class~$\al$; for ease of reading, 
this dependence is dropped from the notation. 

The exact sequence~(\ref{seq_of_Gamma}) of Lemma~\ref{groups_lemma} gives rise 
to an exact sequence of pointed sets (\cite[Proposition 4.5]{mi}):
\[ H^1(B\e, C_\Gamma) \stackrel{\gamma}{\ra} H^1(B\e, N_\Gamma) \stackrel{\delta}{\rightarrow} H^1(B\e, \aut(\Delta)).\]
Consider the class~$\delta(\al)\in H^1(B\e, \aut(\Delta))$; let 
$A$ be the minimal subgroup of~$\aut(\Delta)$ such that 
$\delta(\al)$ is in the image of 
\[H^1(B\e,A)\hookrightarrow H^1(B\e, \aut(\Delta)).\] 
The preimage of~$\delta(\al)$ in~$H^1(B\e,A)$ defines 
a finite \'etale cover~$b\colon {\tilde B}\ra B$. 
Since the cocycle cannot be reduced to a smaller subgroup,
${\tilde B}$ is connected. 

For the rest of this paper, I~fix the type of the simply laced diagram~$\Delta$ 
(equivalently the finite subgroup~$\Gamma$) and the cohomology class~$\alpha$;  
hence also~$A$ and~$\Xi=\Delta/A$ are fixed.  By the exact sequence above, 
$A$ is trivial if and only if~$\alpha$ takes values in~$C_\Gamma$. 

Since~$\Gamma$ is a subgroup of~$\SL(2,\C)$, there is a map 
$\det\colon N_\Gamma\ra\C^*$. Hence~$\al$ gives rise to an induced cocycle
\[\det(\al)\in H^1(B\e, \C^*)\cong\pic(B);\] I~will denote by~$\M$ the  
corresponding line bundle on~$B$. Recall also that 
Proposition~\ref{defs_of_Y}(iii) defines a map $N_\Gamma\ra\GL(\h_\Delta)$.
The image of the cohomology class~$\al\in H^1(B\e, N_\Gamma)$ 
in~$H^1(B\e, \GL(\h_\Delta))\cong H^1(B, \GL(\h_\Delta))$ 
induces a locally free sheaf~$\sH$ of rank~$r$ on~$B$.  

If~$A$ is trivial, this construction gives a class in 
\[H^1(B, \C^*)\subset H^1(B, \GL(\h_\Delta)),\]
and in this case~$\sH\cong\M\otimes\h_\Delta$. A positive root 
$\la\in\Sigma_\Delta^+\subset\h_\Delta$ defines a map 
$\la\colon \h_\Delta\rightarrow\C$ using the inner product, which globalizes to a 
map of bundles 
$\la\colon \sH\rightarrow\M$ over~$B$ and hence to a map on sections
\begin{equation}\label{simple_map_sect}
m_\la\colon H^0(B, \sH)\ra H^0(B,\M).
\end{equation} 

If~$A$ is nontrivial, there is an isomorphism~$b^*\sH \cong b^*\M \otimes \h_\Delta$
on the \'etale cover~$\tilde B$ of~$B$. On the other hand, by Leray
\[ H^0(\tilde B, b^*\sH)\cong H^0(B, b_*b^*\sH)\cong H^0(B, \sH\otimes b_*\O_{\tilde B})
\]
and taking~$A$-invariants, 
\begin{eqnarray*} H^0(B, \sH) & \cong  & H^0(\tilde B, b^*\sH)^A \\ & \cong & H^0(\tilde B, b^*\M)^A\otimes\h_\Delta^A\cong H^0(\tilde B, b^*\M)^A\otimes \h_\Xi.\end{eqnarray*}
Hence if~$\mu\in\Sigma_\Xi^+$ is a positive root in the root system of~$\Xi$, 
there is an induced map 
\begin{equation}\label{map_sect}m_\mu\colon H^0(B, \sH)\ra H^0({\tilde B},b^*\M)^A.
\end{equation}

\subsection{The threefold~$X$}\label{threefold} 

Represent the group cohomology element~$\al\in H^1(B\e, N_\Gamma)$ by a \v Cech 
cocycle~$\{\al_{ij}\in\Gamma(B_{ij}, N_\Gamma)\}$ with respect to an \'etale 
covering~$\{B_l\}$ of~$B$. 
Consider the product
$\C^2\times B_l$, its quotient~$\bar X_l = \C^2/\Gamma\times B_l$ and 
resolution~$X_l= Y\times B_l\ra \bar X_l$ over~$B_l$. By \'etale descent, 
I~can glue the morphisms~$X_l \ra \bar X_l$ over~$B$ 
using the cocyle~$\al$ with values in~$N_\Gamma\subset\aut(Y)$  
(Proposition~\ref{auts_of_ALE}) to get a diagram of quasiprojective varieties 
\[\begin{array}{rcl} X & \stackrel{f}{\longrightarrow}  & \bar X \\
&\mapse{\pi}&\mapdownright{\bar\pi}\\
&&B
\end{array}\]
which up to isomorphism only depends on the cohomology class of~$\al$. 

\begin{propositio}\label{geom_of_X}\ \begin{enumerate} 
\item The singular locus of~$\bar X$ is a section of~$\bar\pi$
and is a curve of cDV singularities of the type of the diagram~$\Delta$. 
The morphism~$f$ is a crepant resolution of singularities. 
The canonical bundle of~$X$ is 
\[ \omega_{X} \cong \pi^*(\omega_B\otimes\M^{-1})\]
and~$H^2(X, \O_X)=0$. 
\item Assume that~$(\Delta, A)\neq (A_{2n}, \Z/2)$ as usual. The exceptional 
locus of~$f$ is a union of irreducible divisors~$D_j$ indexed by the nodes of 
$\Xi$; every such divisor is a smooth geometrically ruled surface 
$\pi_j\colon D_j\ra B_j$. If~$j$ is a node of~$\Xi$ representing 
an~$A$-orbit of size one then~$B_j\cong B$, otherwise~$B_j\cong\tilde B$. 
Two divisors~$D_i, D_j$ intersect if and only if the corresponding nodes~$i,j$ 
of~$\Xi$ are connected by an edge ($m_{ij}>2$). 
\item For every subset~$I$ of the nodes of~$\Xi$, there exists a 
morphism~\[f_I\colon X\ra \bar X_I\] which contracts the divisors
$D_j$ for~$j\in I$ along their rulings and is an isomorphism on 
$X\setminus \cup_{j\in I} D_j$. 
\end{enumerate} 
\end{propositio} 

\begin{figure}[ht]
\centering
\input{conf.pstex_t}
\caption{Dynkin diagrams and configurations of surfaces} 
\end{figure}
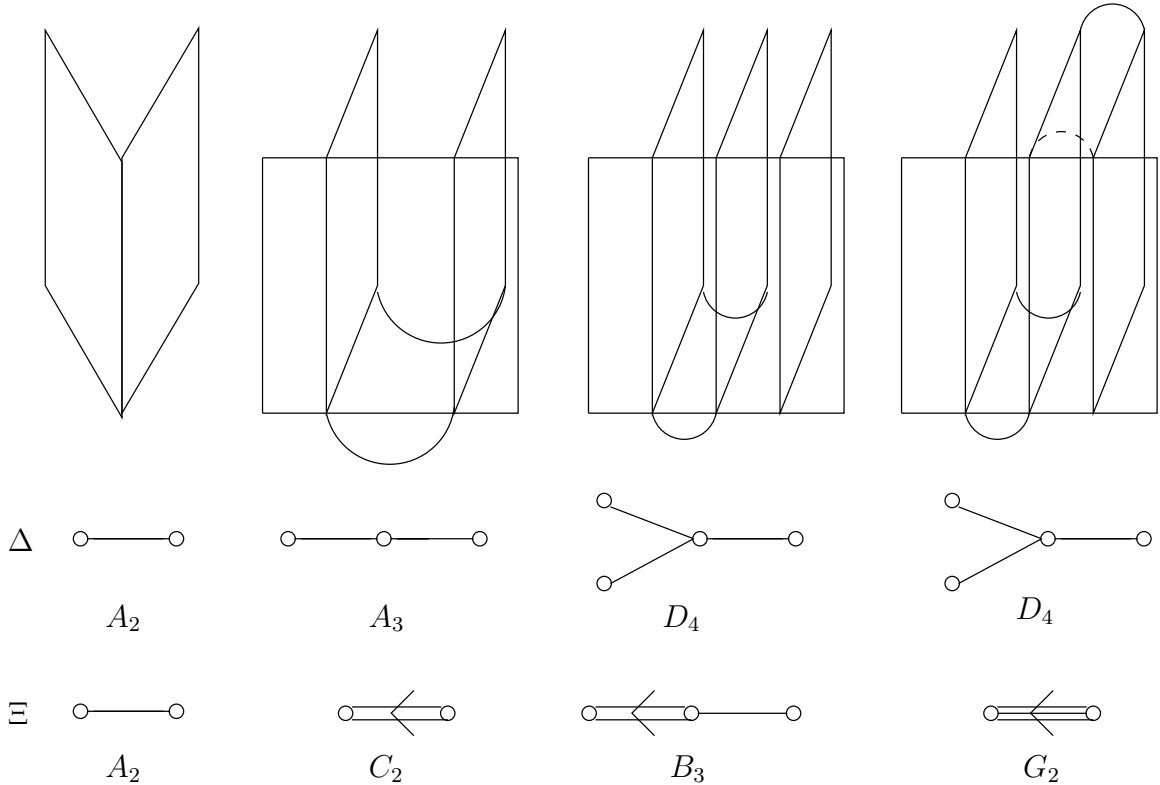

\proof The first statement of (i) is clear. The statement that~$f$ is 
a crepant resolution is local and hence follows from the fact that 
$Y\ra \C^2/\Gamma$ is a crepant resolution. To compute the canonical bundle
of~$X$, note that with~$\pi_j\colon X_j= Y\times B_j\ra B_j$, 
$\omega_{X_j}= \pi_j^*(\omega_{B_j})$, and these line bundles glue together
after a twisting by the inverse of the determinant cocycle. 
$H^2(X,\O_X)=0$ follows from the Leray spectral sequence for~$\pi$ 
and well-known properties of~$Y$. 

For (ii), recall that~$X$ was glued together from \'etale open subsets 
$X_l= Y\times B_l$. If~$\al$ can be represented by a cocycle with 
values in~$C_\Gamma$, in other words if~$A$ is trivial, 
then the glueing process will not permute the exceptional divisors~$\{E_j\}$ 
in~$Y$ by Proposition~\ref{auts_of_ALE}. 
Hence~$\{E_j\times B_l\}$ will glue for every~$j$ 
to a smooth exceptional divisor~$D_j$ ruled over the curve~$B$, and these 
surfaces will intersect as dictated by the diagram~$\Xi=\Delta$.
If~$A$ is nontrivial, then the set of exceptional lines~$\{E_j\}$ is acted on
by monodromy over~$B$, this action being given by the cocyle~$\delta(\al)$
with values in~$A$. The lines~$E_i$ corresponding to nodes fixed by the action of 
$ A$ can still be glued globally over~$B$, leading to exceptional 
divisors in~$X$ ruled over~$B$. However, nontrivial orbits
$\{E_{j_i}\}$ of exceptional curves under the 
$ A$-action are glued together to an irreducible exceptional divisor, 
where the glueing is governed by the cohomology class~$\delta(\al)$; hence
the corresponding surfaces are ruled over the \'etale cover~${\tilde B}$ of~$B$. 
Finally the morphism~$f_I\colon X\ra \bar X_I$ can be glued over~$B$ from the morphism
$G_J\colon Y\ra \bar Y_J$, contracting the exceptional curves on~$Y$ given by the 
$A$-invariant set~$J$ of nodes of~$\Delta$ lying over~$I$. 
\ethrm

\begin{remar}\rm\label{specialcase} In the special case
$(\Delta, A)=(A_{2n}, \Z/2)$ 
the exceptional divisors~$D_i$ are still indexed by nodes of the quotient diagram
$\Delta/A$, defined to be the~$A_n$-diagram with a marked 
node at one end corresponding to the adjacent~$\aut(\Delta)$-orbit of nodes. 
However, the marked node~$n$ of~$\Xi$ corresponds to a singular 
exceptional surface. It is an irreducible non-normal surface 
$\pi_n\colon D_n\ra B$ whose double locus is a section and whose fiber over 
any point $b\in B$ is a line pair. The main results of this
paper do not apply in this special case; 
see~\cite[Remark 4.5]{eg} for further discussion. 
\end{remar} 

\subsection{A family of deformations of~$X$}\label{sec!defs_of_X}

Representing the cohomology class~$\al$ as a \v Cech cocyle again with respect
to an \'etale covering~$\{B_l\}$ of~$B$, 
there are isomorphisms~$\sH|_{B_l}\cong \h_\Delta\times B_l$. 
Using these isomorphisms, pull back 
the universal deformation space~$\Y\rightarrow\h_\Delta$ 
of~$Y$ from Proposition~\ref{defs_of_Y} to a family of surfaces
$\tilde\X_l\rightarrow\sH|_{B_l}$. Glue these families 
over~$B$ using the identification given by the cocyle with values in 
$N_\Gamma<\aut(\Y)$ (compare Proposition~\ref{defs_of_Y}(iii))
to get a global family of surfaces 
$\tilde\X\rightarrow\sH$ over the total space of the vector 
bundle~$\sH$. Finally use the tautological map 
$B\times H^0(B, \sH)\rightarrow \sH$ to pull back~$\tilde\X$ to 
a family~$\X$ over~$B\times H^0(B, \sH)$. This leads to a diagram 
\[ \begin{array}{ccc}\X & \rightarrow &\tilde\X\\
\downarrow && \downarrow \\B\times T& \rightarrow & \sH \\
\downarrow \\ T \end{array}\]
with $T=H^0(B, \sH)$. 
The composite of the left-hand vertical maps gives rise to a morphism
$e\colon\X\rightarrow T$, which is a smooth family of
threefolds by construction. 
Using the projection to~$B$ shows that for~$s\in T$, the fiber 
$X_s=e^{-1}(s)$ admits a smooth map~$\pi_s\colon X_s\ra B$. 
The central fiber~$e^{-1}(0)$ corresponds to the zero section of
the bundle~$\sH$; it is obtained by glueing varieties~$Y\times B_l$, 
coming from the central fiber of~$\Y\ra \h_\Delta$, over the curve~$B$
as dictated by the cohomology class~$\al$. Thus~$e^{-1}(0)\cong X$. 
Note that at this point~$T$ may well be infinite dimensional, but the 
meaning of the following statements should be obvious also in this case. 

\begin{propositio}\label{defs_of_X} The family~$e\colon \X\ra T$ with central 
fiber~$X\cong e^{-1}(0)$ has the following properties:
\begin{enumerate} 
\item The Kodaira--Spencer map of the family is injective at~$0\in T$. 
\item For any set~$I$ of nodes of~$\Xi$, there is a contraction morphism 
$F_I\colon \X\ra\bar\X_I$ over~$T$ with central fiber~$f_I\colon X\ra \bar X_I$.
\item The group~$W_\Xi$ acts on the base~$T$ with the following 
properties: 
\begin{enumerate} 
\item for~$i$ a node of the Dynkin diagram~$\Xi$, the fixed 
locus~$T_i={\rm Fix}(\rho_i)$ of the reflection 
$\rho_i\in W_\Xi$ acting on~$T$ is exactly the locus of points 
$s\in T$ for which the map~$f_{i,s}\colon X_s\ra X_{i,s}$ 
(the fiber of~$F_i$ at~$s\in T$) contracts a surface;   
\item for a reflection~$w\in W_\Xi$, if~$I$ is the set of simple roots 
reflected by~$w$, there exists a diagram 
\[\begin{array}{ccc}
\X &\stackrel{\theta_w}\dashrightarrow & \X \\
\Big\downarrow && \Big\downarrow \\
\bar \X_I & \stackrel{\sim}\longrightarrow & \bar \X_I \\
\Big\downarrow && \Big\downarrow \\
T &\stackrel{w}\longrightarrow & T. 
\end{array}\]
The relative birational morphism~$\theta_w\colon \X\dashrightarrow\X$ over~$T$ restricts 
to a birational morphism~$\theta_{w,s}\colon X_s\dashrightarrow X_{w(s)}$ on fibers.
\end{enumerate} 
\end{enumerate} 
\end{propositio}
\proof The Kodaira--Spencer map of the family~$e\colon\X\ra T$ at~$0\in T$ is a map
\[ \phi\colon \Theta_{T,0} \cong H^0(B,\sH)\longrightarrow H^1(X, \Theta_X)
\] 
where~$\Theta_{T,0}$ is the holomorphic tangent space of~$T$ at~$0\in T$, and
$\Theta_X$ is the holomorphic tangent bundle of~$X$. This map sits in a 
composition 
\begin{equation}\label{ksrow}
H^0(B,\sH)\rightarrow H^1(X, \Theta_X) \rightarrow H^0(B, R^1\pi_*\Theta_X)\rightarrow H^0(B,  R^1\pi_*\Theta_{X/B})
\end{equation}
where the second map comes from the Leray spectral sequence, and the last map 
comes from the natural map of sheaves~$\Theta_X\rightarrow\Theta_{X/B}$ on~$X$. 
On the other hand, it is easy to check from the construction that 
the Kodaira--Spencer map~$\h_\Delta\stackrel{\sim}\rightarrow H^1(Y, \Theta_Y)$ for~$Y$ 
globalizes to an 
isomorphism between the sheaf~$\sH$ and the sheaf~$R^1\pi_*\Theta_{X/B}$ on~$B$, 
and the composite~$H^0(B,\sH) \rightarrow H^0(B,  R^1\pi_*\Theta_{X/B})$
of the maps in~(\ref{ksrow}) 
is the induced isomorphism. Hence the Kodaira--Spencer map of~$e\colon\X\ra T$
is injective. 

For (ii), let~$J$ be the subset of exceptional curves in~$Y$ corresponding 
to the set of nodes~$I$ of~$\Xi$, and let~$G_J\colon \Y\ra \Y_J$ be the corresponding 
contraction from Proposition~\ref{defs_of_Y}(vi). 
$J$ is fixed under the monodromy action by~$A$, 
so I~can glue~$G_J$ over~$B$ according to the cocycle~$\al$ to get a morphism 
$F_I\colon \X\ra \bar \X_I$ with central fiber~$f_I$. 

For (iii), note that~$W_\Xi<W_\Delta$ acts on~$\h_\Delta$ and by Remark~\ref{acts}
on p.~\pageref{acts} this action commutes with that of~$A$ 
and of course the scalars. Hence it acts on the vector bundle~$\sH$ over~$B$ 
(trivial action on~$B$) and so on~$T=H^0(B, \sH)$. 

To show property (a) of this action, assume first that~$A$ is 
trivial. Fix a section~$s\in T$ of~$\sH$ and a node~$i$ of~$\Xi=\Delta$. 
Note that by Proposition~\ref{defs_of_Y}(v), the surface~$X_{s,p}=\pi_s^{-1}(p)$ 
over~$p\in B$ contains a deformation of the rational curve 
$E_i\subset Y$ if and only if~$m_{r_i}(s)$ vanishes at~$p\in B$, where
$m_{r_i}$ is the map in~(\ref{simple_map_sect}). 
The contraction~$f_{i,s}$ contracts a surface if and only if this happens
at every point~$p\in B$, in other words if~$m_{r_i}(s)=0$. However, this simply
says that in every fiber~$s(p)$ is on the reflection hyperplane~$\Pi_{r_i}\subset\h_\Delta$
or equivalently that it is fixed by~$r_i$. Hence the fixed locus of~$r_i$ 
acting on~$T$ is exactly the locus where~$f_{i,s}$ contracts a surface. 

If~$A$ is nontrivial, the node~$i$ of~$\Xi$ corresponds to an~$A$-orbit 
$\{i_j\}$ of nodes of~$\Delta$. For~$s\in T$, the surface fiber~$X_{s,p}$ over 
$p\in B$ will contain deformations of the exceptional curves~$E_{i_j}$ of~$Y$ 
if and only if~$m_{\rho_i}(s)$ vanishes at the preimages 
$b^{-1}(p)\subset\tilde B$ of~$p\in B$, where~$m_{\rho_i}$ is the map 
in~\ref{map_sect}. Hence~$f_{i,s}$ contracts a surface if and only if 
$m_{\rho_i}(s)=0$; equivalently, if in every fiber~$s(p)\in\sH_p\cong\h_\Delta$ 
lies on all reflection hyperplanes~$\Pi_{r_{i_j}}\subset\h_\Delta$. 
This however happens if and only if the section~$s$ is fixed by 
$\rho_i=\prod r_{i_j}\in W_\Xi\subset W_\Delta$. 

The isomorphism~$\bar\X_I\stackrel{\sim}\ra\bar\X_I$ fitting into the diagram of (b)
for a reflection~$w\in W_\Xi$ is the pullback of the diagram of 
Proposition~\ref{defs_of_Y}(iv) and it naturally induces a birational map between
resolutions; the details are left to the reader. 
\ethrm 

\begin{remar}\rm  The~$W_\Xi$-action on~$T$ will be essential in what
follows. I~continue calling an element~$w\in W_\Xi$ a reflection
if it is a reflection on~$\h_\Xi$, that is if it corresponds to a positive 
root in~$\Sigma_\Xi^+$. Of course in general the fixed locus of~$w$ on~$T$ is 
not a hyperplane. 
\end{remar}

If the linear system~$\M$ on~$B$ is small, then the family 
$e\colon \X\ra T$ can be rather uninteresting (for example,~$T$ could be a point). 
Under an extra assumption, a lot more geometry emerges. 

\begin{propositio} \label{many_defs_of_X}
Assume that~$\M$ is a moving linear system on~$B$.
\begin{enumerate}
\item For a general point~$s\in T$, the exceptional locus of 
$f_s\colon X_s\ra \bar X_s$ consists of a finite number of disjoint smooth rational 
curves, naturally indexed by positive roots 
$\mu\in\Sigma_\Xi^+$, with normal bundle $\O_{\P^1}(-1,-1)$.
\item For such general~$s\in T$ and a reflection~$w\in W_\Xi$, the 
birational map~$\theta_{w,s}\colon X_s \dashrightarrow X_{w(s)}$ is a flop, 
flopping exactly those~$(-1,-1)$-curves which are indexed 
by positive roots mapped to negative roots by~$w$.
\item For~$w=\rho_i$ with~$i$ a node of the Dynkin diagram~$\Xi$, 
there are two possibilities: either~$s\in T_i$ and 
$\theta_{\rho_i,s}\colon X_s \dashrightarrow X_{\rho_i(s)}$ is the 
identity isomorphism, or~$s\in T\setminus T_i$ and the birational map 
$\theta_{\rho_i,s}\colon X_s \dashrightarrow X_{\rho_i(s)}$ flops a
disjoint union of rational curves indexed by the simple root~$\mu_i$, with normal 
bundle $\O_{\P^1}(-1,-1)$ or $\O_{\P^1}(0,-2)$. 
\end{enumerate} 
\end{propositio}
\proof To prove (i), assume first that~$A$ is trivial. 
Let~$p\in B$ a closed point and~$X_{s,p}=\pi_s^{-1}(p)$ the corresponding 
surface. According to Proposition~\ref{defs_of_Y}, the surface~$X_{s,p}$ is affine
if and only if~$s(p)\in\h_\Delta=\sH_s$ does not lie on any reflection hyperplane
$\Pi_w$, and it contains a single smooth rational curve if it lies in a unique 
such hyperplane. Said invariantly, using the map~(\ref{simple_map_sect}), 
$X_{s,p}$ is affine if and only if~$m_\la(s)$ does not vanish at~$p$ for any 
$\la\in \Sigma_\Delta^+$, and it contains
a unique curve if~$m_\la(s)$ vanishes at~$p$ for a unique 
$\la\in \Sigma_\Delta^+$. 

I~claim that for general~$s\in T$, there is a finite set of~$p\in B$ such that 
$m_\la(s)$ vanishes at~$p$ for some positive root~$\la$, and at such a~$p$ there 
is a unique $\la\in \Sigma_\Delta^+$ 
with~$m_\la(s)(p)=0$. This shows that for such~$s$ there is a finite 
number of disjoint rational curves in~$X_s$, which are indexed by positive roots. 

To show the claim, take two sections 
$t, t'\in H^0(B, \M)$ which have a disjoint set of simple zeroes; 
as~$|\M|$ has no base points, this is possible. Choose also 
$h, h'\in\h_\Delta$ so that for~$\la, \la'\in\Sigma_\Delta^+$ 
different positive roots, 
\begin{equation}\label{detcond}\langle\la, h\rangle\langle\la', h'\rangle-\langle\la, h'\rangle\langle\la', h\rangle\neq 0.\end{equation} 
Setting~$s=h\otimes t + h'\otimes t'$, the section 
$m_\la(s)$ of~$\M$ 
vanishes only at finitely many points of~$B$ for any positive root~$\la$. Also,
if~$m_\la(s)$ vanishes at the same point 
$p\in B$ for two different roots, then by condition~(\ref{detcond}), 
I~get that~$t(p)=t'(p)=0$ which contradicts the choice of~$t, t'$. 

The proof in the general case, when~$A$ is nontrivial, is similar. 
In this case, the claim is that for general~$s\in T$ there is a finite number of 
points~$p\in B$ such that~$m_\la(s)\in H^0({\tilde B},b^*\M)$ 
vanishes at a point~$q\in {\tilde B}$ lying over~$p\in B$ and for a positive root
$\mu\in\Sigma_\Xi^+$, and at such points~$p$, vanishing happens for a unique~$\mu$. 
The claim follows by considering
\[b^*t\otimes h+ b^*t'\otimes h'\in H^0({\tilde B}, b^*\M)^A\otimes\h_\Xi\cong H^0(B, \sH)\]
for sufficiently general~$h,h'\in\h_\Xi$ and~$t,t'\in H^0(B,\M)$
with no common zeros. 

By construction and the discussion above, a small analytic neighbourhood 
of every rational curve on the general fiber~$X_s$ looks like the standard 
one-dimensional deformation of a~$(-2)$-curve in a surface, the deformation 
direction being transversal to the hyperplane along which the curve 
deforms. In other words, locally near the curve, the threefold looks like 
a small resolution of the ordinary threefold double point $\{xy=z^2-t^2\}$. 
It is well known that the normal bundle of such a curve is~$\O_{\P^1}(-1,-1)$. 
This also proves statement (ii): by Proposition~\ref{defs_of_X}(iii), the 
birational map~$\theta_{w,s}$ factors 
as~$X_s\rightarrow \bar X_{I, s}\leftarrow X_{w(s)}$, where the first map 
contracts exactly those~$(-1,-1)$-curves which are indexed by positive roots
mapped to negative ones by~$w$; locally analytically the two maps give the two 
small resolutions of the resulting ordinary double points, in other words the flop. 

To prove (iii), note that one possibility is~$m_{\rho_i}(s)=0$, 
hence~$s\in T_i$ and the birational map $\theta_{\rho_i,s}$ is the identity. 
Otherwise~$m_{\rho_i}(s)$ vanishes at a finite set of points $\{p_l\}\subset B$, 
and the map~$X_s\ra \bar X_{i,s}$ contracts a finite number of disjoint rational
curves contained in the surfaces~$X_{s,p_l}$. If at $p_l$ the section 
$m_{\rho_i}(s)$ meets the zero-section of $\M$ transversally, 
then the singularity on $\bar X_{i,s}$ is analytically isomorphic 
to $\{xy=z^2-t^2\}$ and the corresponding rational 
curve has normal bundle $\O_{\P^1}(-1,-1)$. If the intersection of the
two sections is tangential of order $n>1$ then the singularity 
on~$\bar X_{i,s}$ is analytically isomorphic to $\{xy=z^2-t^{2n}\}$ 
and the normal bundle is $\O_{\P^1}(0,-2)$. 
In any case,~$\theta_{\rho_i,s}$ flops all these curves (compare~\cite{pagoda}).  
\ethrm

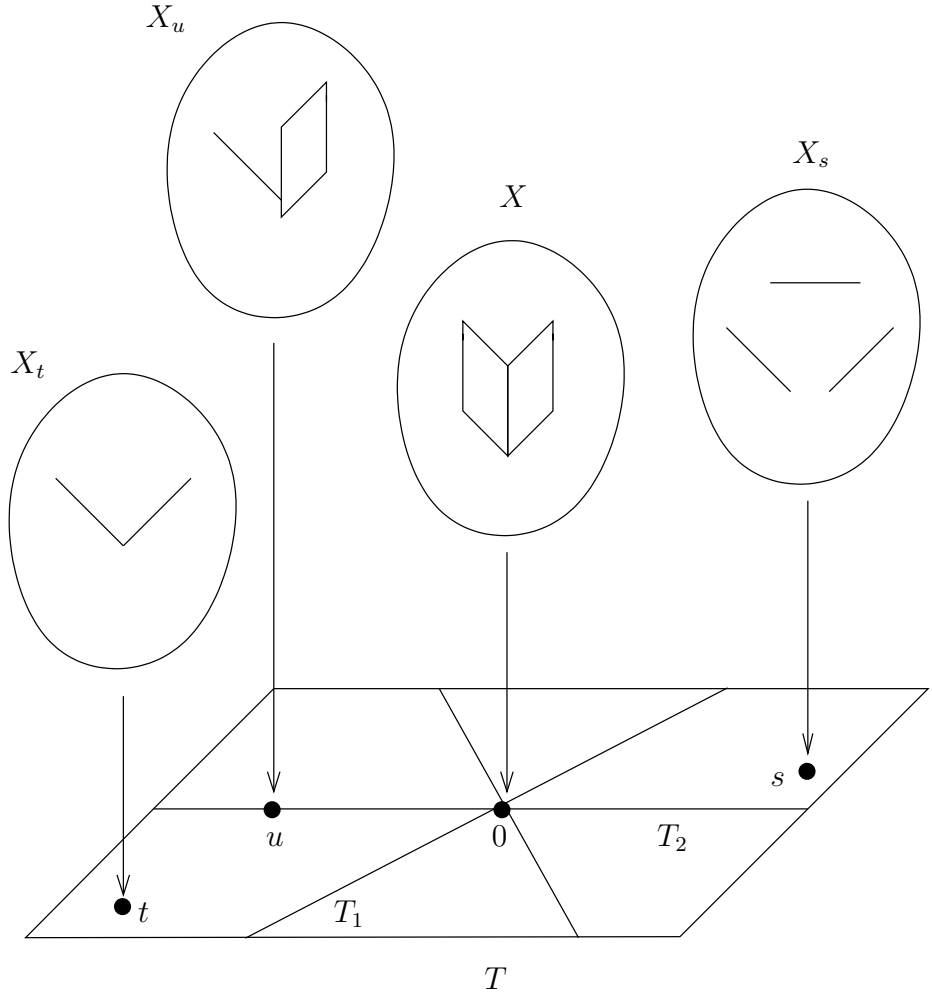
\begin{figure}[ht]
\centering
\input{loci.pstex_t}
\caption{Possible exceptional loci for type~$A_2$}
\end{figure}

\begin{remar} \rm 
Points~$s\in T$ in the base with~$f_s$ having an exceptional locus
consisting of a disjoint union of smooth~$(-1,-1)$-curves as in (i) will be called 
{\em sufficiently general}; on Figure~2,~$s\in T$ is a sufficiently 
general point but~$t\in T$ is not. 
As the proof shows, the locus of sufficiently
general points is a non-empty Zariski open subset of~$T$. 

Note also that it can perfectly well happen that~$\M$ is the trivial line
bundle on~$B$ and its only sections are the constants.
All the statements of the above discussion remain true, with the
small proviso that the maps~$\theta_{w,s}$ flop an empty set of curves, 
in other words they are isomorphisms for all~$s$. This phenomenon
(in the projective case) is well known in the literature; the nontrivial 
birational contraction~$f\colon X\rightarrow \bar X_i$ on the central fiber 
deforms to an isomorphism
$f_s\colon X_s\stackrel{\sim}\ra \bar X_{i,s}$, hence the cone of ample divisors 
jumps in the family. The variety~$X$ is Calabi--Yau if and only if 
$\M\cong\omega_B$ by Proposition~\ref{geom_of_X}(i), 
hence such surfaces are elliptic ruled surfaces in 
Calabi--Yau threefolds. Compare~\cite{wi}--\cite{wi_ell}.
\end{remar}

\section{Derived categories and equivalences in families}\label{sec!dq}

\subsection{Kernels and Fourier--Mukai functors}

If~$X$ is a smooth quasiprojective variety, let~$\db{X}$ denote 
the bound\-ed derived category of coherent sheaves on~$X$. 
A {\it kernel} (derived correspondence) between smooth quasiprojective varieties 
$X_i$~$(i=1,2)$ is an object~$U\in\db{X_1\times X_2}$, whose
support is proper over both factors. 
There is a composition product on kernels given for  
$U\in\db{X_1\times X_2}$ and~$V\in\db{X_2\times X_3}$ by the standard formula
\[ U\circ V= \push p_{13*} \left(p_{23}^*(V)\tensor p_{12}^*(U)\right)\in \db{X_1\times X_3};
\]
here~$p_{ij}\colon X_1\times X_2\times X_3\rightarrow X_i\times X_j$ 
are the projection maps and the pullbacks are ordinary pullbacks since 
$p_{ij}$ is flat. A kernel~$U\in\db{X_1\times X_2}$ 
is {\it invertible}, if there is a kernel~$V\in\db{X_2\times X_1}$ 
such that the products~$U\circ V$ and~$V\circ U$ are 
isomorphic in~$\db{X_i\times X_i}$ to~$\O_{\Delta_{X_i}}$, the (complexes
consisting of) the structure sheaves of the diagonals. 
A kernel~$U\in\db{X_1\times X_2}$ defines a functor
\[\Psi^U\colon \db{X_2}\rightarrow\db{X_1}
\]
by 
\[\Psi^U(-)= \push p_{1*}(U\tensor p_2^*(-)), 
\]
with~$p_i\colon X_1\times X_2\ra X_i$ the projections. If~$U$ is invertible then 
$\Psi^U$ is a {\it Fourier--Mukai functor}, an equivalence of triangulated 
categories. Let $\daut{X}$ denote the group of invertible kernels on $X$. 

\subsection{Kernels in families}\label{dcats} 

Suppose that~$\pi_i\colon \X_i\ra S_i$,~$i=1,2$ are smooth families. 
A {\em relative kernel for~$(\pi_1,\pi_2)$} is a pair~$(U,\phi)$, where  
\begin{itemize}
\item $\phi\colon S_1\rightarrow S_2$ is a isomorphism, 
giving rise to the fiber product diagram 
\[\begin{array}{ccc}
\X_1\times_\phi \X_2 & \stackrel{p_2}{\longrightarrow} & \X_2 \\
\mapdownright{p_1} &\stackrel{\pi_{12}}{\searrow}& \mapdownright{\phi^{-1}\circ\pi_2} \\
\X_1 & \stackrel{\pi_1}{\longrightarrow} &  S_1
\end{array}\]
and
\item $U\in\db{\X_1\times_\phi \X_2}$, whose derived restriction to the 
fiber of~$\pi_{12}$ over every~$s\in S_1$ is isomorphic to an object in 
$\db{X_{1,s}\times X_{2,\phi(s)}}$ with proper support over both factors.  
\end{itemize}

If~$\pi_i\colon \X_i\ra S_i$,~$i=1,2,3$ are smooth families and~$(U,\phi), (V,\psi)$ 
relative kernels for~$(\pi_1,\pi_2), (\pi_2,\pi_3)$, 
then there is a composition on kernels defined
as follows. Let~$\chi=\psi\circ\phi\colon S_1\ra S_3$, and let 
\begin{equation}\label{twice}\pi_{123}\colon\tilde\X= (\X_1\times_\phi\X_2)\times_\chi\X_3\rightarrow S_1\end{equation}
be the twice-fiber product. 
There are maps~$p_{12}\colon \tilde\X\ra\X_1\times_\phi\X_2$
and~$p_{13}\colon \tilde\X\ra\X_1\times_\chi\X_3$ which commute with the maps to~$S_1$. 
On the other hand, it is easy to check that 
$(\X_1\times_\phi\X_2)\times_\chi\X_3\cong \X_1\times_\phi(\X_2\times_\psi\X_3)$ 
hence there is a map~$p_{23}\colon \tilde\X\ra\X_2\times_\psi\X_3$ which satisfies
$\pi_{23}\circ p_{23}= \phi\circ\pi_{123}$. Hence finally I~can put
\[W= \push p_{13*} (p_{23}^*(V)\tensor p_{12}^*(U))\]
and set 
\[ (V,\psi)\circ (U,\phi) = (W, \chi).\]
It is easy to see that~$(W, \chi)$ is a relative kernel for~$(\pi_1, \pi_3)$. 

For reference I~record the composition of three kernels, leaving the obvious 
generalization to the reader. Let~$\pi_i\colon\X_i\ra S_i$
be families for~$i=1,\ldots, 4$, and~$(U,\phi), (V,\psi), (T,\eta)$ relative 
kernels for~$(\pi_1, \pi_2)$, $(\pi_2, \pi_3)$  and $(\pi_3, \pi_4)$. Let
\[\widehat\X= \X_1\times_\phi(\X_2\times_\psi(\X_3\times_\eta\X_4))\]
with maps~$s_{12}\colon \widehat\X\ra \X_1\times_\phi\X_2, s_{23}\colon \widehat\X\ra \X_2\times_\psi\X_3, s_{34}\colon \widehat\X\ra \X_3\times_\eta\X_4$ and~$s_{14}\colon \widehat\X\ra \X_1\times_{\eta\circ\psi\circ\phi}\X_4$. Then 
\begin{lemm} \[ (T,\eta)\circ (V,\psi)\circ (U,\phi) \cong (\tilde U, \eta\circ\psi\circ\phi),\] 
where
\[\tilde U = \push s_{14*} (s_{34}^*(T)\tensor s_{23}^*(V)\tensor s_{12}^*(U)).\]
\label{three_ker}\end{lemm} 
\ethrm

A relative kernel~$(U,\phi)$ for~$(\pi_1, \pi_2)$ is called {\it invertible},
if there is a relative kernel~$(V,\phi^{-1})$ for~$(\pi_2, \pi_1)$, with the
property that 
the compositions~$(U,\phi)\circ(V,\phi^{-1})$ and~$(V,\phi^{-1})\circ(U,\phi)$
are isomorphic to the relative kernels $(\O_{\Delta_{\X_i}}, {\rm id})$. 

Let $U_s\in\db{ X_{1,s}\times X_{2,\phi(s)}}$ denote the derived restriction of 
the relative kernel $U$ to fibers of $\pi_{12}$.
\begin{propositio} If a relative kernel~$U\in\db{\X_1\times_\phi\X_2}$ 
for~$(\pi_1, \pi_2)$ is invertible then the restricted kernel
$U_s\in\db{ X_{1,s}\times X_{2,\phi(s)}}$ is invertible for every~$s\in S_1$. 
Conversely if the restricted kernel is invertible for every~$s\in S_1$, then 
every~$s\in S_1$ has a neighbourhood~$s\in T\subset S_1$ such that the
relative kernel restricted to the families~$\pi_1^{-1}(T)\ra T$ and 
$\pi_2^{-1}(\phi(T))\ra\phi(T)$ is invertible. 
\label{fiberwise}
\end{propositio}
\proof  If~$x_{i,s}\colon X_{i,s}\hookrightarrow\X_i$ denotes the inclusion, 
then 
\[\pull x_{i,s}^*\O_{\Delta_{\X_i}} \cong \O_{\Delta_{X_{i,s}}}\]
as the maps~$\pi_i$ are flat. Hence if~$U$ is invertible with inverse~$V$, 
then~$U_s$ is invertible with inverse~$V_s$. 

Conversely, take a relative kernel~$U\in\db{\X_1\times_\phi\X_2}$ and suppose that
the restrictions are all invertible. Let 
\[ V=\rhom(U, \O_{\X_1\times_\phi\X_2})\otimes p_1^*\omega_{\X_1/S_1}[n]\in \db{\X_2\times_\phi\X_1}
\]
where~$n$ is the dimension of the fibers~$X_{i,s}$. 
Then by standard adjunctions the functor~$\Psi^{(V,\phi^{-1})}$ is 
right adjoint to the functor~$\Psi^{(U,\phi)}$. 
Hence~$\Psi^{V_s}$ is right adjoint to 
$\Psi^{U_s}$ on the fibers. However, adjoints are unique, so~$V_s$ is the 
inverse of the kernel~$U_s$. In other words, 
\[\pull z_{1,s}^*\left(U\circ V\right) \cong U_s\circ V_s \cong \O_{\Delta_{X_{1,s}}}\in \db{X_{1,s}\times X_{1,s}},\]
where~$z_{1,s}\colon X_{1,s}\times X_{1,s} \rightarrow \X_1\times_S \X_1$ is the 
inclusion. By~\cite[Lemma 4.3]{bridgeland}, 
this implies that~$U\circ V$ is a 
sheaf on~$\X_1\times_S \X_1$, flat over~$S$. 
Take~$s\in S_1$, then the natural map 
\[ H^0(\X_1\times_S\X_1, U\circ V)\ra H^0(X_{1,s}\times X_{1,s}, U_s\circ V_s)\]
is surjective, hence there is a map of sheaves
\[\O_{\X_1\times_S\X_1}\ra U\circ V\]
which is surjective at~$s\in S_1$. So this map is surjective over a neighbourhood
$T$ of~$s\in S_1$ and restricted to that neighbourhood,~$U\circ V$
is a structure sheaf of a subscheme, fiberwise isomorphic to the diagonal 
$\Delta_{X_{1,s}}\subset X_{1,s}\times X_{1,s}$. By Lemma~\ref{isom_for_subshemes}
below,~$U\circ V$ restricted over~$T$ 
is isomorphic to the structure sheaf of the relative diagonal in the fiber product
$\X_1\times_S\X_1$. To conclude, repeat the argument with~$V\circ U$ and 
take the intersection of the resulting open sets.  
\ethrm

Let now~$\pi\colon \X\rightarrow S$ be a fixed smooth family, and consider relative 
kernels~$(U, \phi)$ where~$\phi\colon S\ra S$ is an automorphism of the base. 
Let~$M(\X/S)$ be the set of such pairs up to isomorphism. 
$M(\X/S)$ has a monoid multiplication
given by composition with a two-sided unit~$({\rm id}_S, \O_{\Delta_\X})$. 
Let~$\daut{\X/S}$ be the group of invertible elements of the monoid~$M(\X/S)$, 
the {\em group of relative equivalences} of the family~$\X\ra S$. 

\subsection{Three lemmas} 

I~record some auxiliary results on sheaves and kernels. 

\begin{lemm} \label{arrows_from_O} Suppose that~$\pi\colon\X\ra S$ is a 
smooth family,\[(U,\phi), (V, \psi)\in M(\X, S)\] are relative kernels with 
composite~$(W, \chi)$ and  
assume further that the sheaf $V$ is isomorphic to the structure sheaf 
$i_*\O_Y$ of a subscheme 
$i\colon \Y\hookrightarrow \X\times_\psi\X$. Then every map 
$\O_{\X\times_\phi\X}\ra U$ in~$\db{\X\times_\phi\X}$ induces a map
$\O_{\X\times_\chi\X}\ra W$ in ~$\db{\X\times_\chi\X}$. 
\end{lemm} 
\proof Since~$V$ is the structure sheaf of a subscheme, 
$p_{23}^*(V)\cong j_*\O_{\tilde\Y}$ 
is also the structure sheaf of the 
subscheme~$j\colon \tilde\Y\hookrightarrow \tilde\X$, where $\tilde\X$ is the
fiber product of~(\ref{twice}).
But by the projection formula
\[p_{23}^*(V)\tensor p_{12}^*(U)\cong j_*\O_{\tilde\Y}\tensor p_{12}^*(U)\cong j_*(\pull j^*p_{12}^*(U)).\]
On the other hand, by functoriality, a map~$\O_{\X\times_\phi\X}\ra U$
induces a map~$\O_{\tilde\X}\ra p_{12}^*(U)$, hence a map 
$\O_{\tilde\Y}\ra\pull j^*p_{12}^*(U)$ and hence a map 
$j_*(\O_{\tilde\Y})\ra j_*(\pull j^*p_{12}^*(U))$. Composing this with 
the natural map~$\O_{\tilde\X}\ra j_*(\O_{\tilde\Y})$ I~get a map
$\O_{\tilde\X}\ra p_{23}^*(V)\tensor p_{12}^*(U)$. But~$\O_{\tilde\X}\cong p_{13}^*(\O_{\X\times_\chi\X})$, so by adjointness I~finally get a map
\[\O_{\X\times_\chi\X}\ra \push p_{13*}(p_{23}^*(V)\tensor p_{12}^*(U))=W.\] 
\ethrm 

\begin{lemm} Let~$f\colon X\rightarrow Y$ be a projective morphism of 
quasiprojective varieties. Let
$i\colon Z\hookrightarrow X$ be a reduced subscheme of~$X$ and 
$j\colon T\hookrightarrow Y$ its reduced image under~$f$. Assume that all fibers of
$f|_Z\colon Z\rightarrow T$ are projective spaces. 
Then~$\push f_*(i_*\O_Z)\cong j_*\O_T$. 
\label{pushdown}\end{lemm} 
\proof As~$i$ and~$j$ are closed immersions, 
the Grothendieck spectral sequences
for~$\push (f\circ i)_*$ and~$\push (j\circ f|_Z)_*$ degenerate, hence 
\[\push f_*(i_*\O_Z)\cong \push (f\circ i)_*(\O_Z)\cong \push (j\circ f|_Z)_*(\O_Z)\cong j_*(\push f|_{Z*}\O_Z).\]
On the other hand,~$f|_{Z*}\O_Z\cong\O_T$ since fibers of~$f|_Z$ are connected, 
and~$R^if|_{Z*}\O_Z =0$ for~$i>0$ by the Theorem on Formal Functions and the fact 
that the higher cohomologies of~$\O$ on projective spaces vanish. 
This proves the statement. 
\ethrm 

\begin{lemm} \label{isom_for_subshemes} Let~$e\colon \X\rightarrow S$ be a smooth 
family. Assume that~$U_1$,~$U_2$ are sheaves on~$\X$, flat over~$S$, with 
surjective maps~$g_i\colon \O_X\rightarrow U_i$ for~$i=1,2$ 
(i.e.~the~$U_i$ are structure sheaves of subschemes). 
Suppose further that for some dense open~$S^0\subset S$, if 
$s\in S^0$ then there is an isomorphism~$U_{1,s}\cong U_{2,s}$ which is compatible
with the maps~$g_{i,s}\colon \O_{X_s}\rightarrow U_{i,s}$. 
Then~$U_1\cong U_2$ as sheaves on~$\X$, compatibly with the maps~$g_i$. 
\end{lemm} 
\proof For~$i=1,2$ the structure sheaves~$U_i$ of subschemes of~$\X$ give rise to 
morphisms 
\[ \phi_i\colon S\ra {\rm Hilb}(\X/S)
\]
over~$S$, where~${\rm Hilb}(\X/S)\ra S$ represents the Hilbert functor 
of the quasiprojective morphism~$\X\ra S$ (this is constructed using 
a projective completion~$\X\hookrightarrow \bar \X\ra S$ along the fibers), 
such that the morphisms~$g_i$ are pullbacks of a universal surjection
$\O_{\X\times_S{\rm Hilb}(\X/S)}\ra \U$. By the condition on restrictions, 
$\phi_1|_{S^0}=\phi_2|_{S^0}$. But the Hilbert scheme is separated, so the maps 
$\phi_1$ and~$\phi_2$ coincide. Hence~$U_1\cong U_2$ compatibly with the maps
$g_i$, since they are pullbacks of~$\U$ along the same map. 
\ethrm 

\section{Artin group actions on derived categories}\label{sec!main}

\subsection{Relative equivalences for threefolds containing ruled surfaces}
\label{sub!rel} 

Recall the family of threefolds~$e\colon \X\ra T$ constructed in Section~\ref{threefold}
together with the action of the reflection group~$W_\Xi$ on the base~$T$. 
Note that~$T$ can a priori be infinite dimensional. For any finite dimensional
$W_\Xi$-invariant vector subspace~$Q\subset T$, I~can consider the restricted 
family~$e_Q\colon \X_Q\rightarrow Q$ which I~will simply denote by~$e\colon \X\ra Q$.
The central fiber of this family is still~$e^{-1}(0)\cong X$. 
Restrict all contractions and the~$W_\Xi$-action to~$Q$. 

For every node~$i$ of the diagram~$\Xi$ there is a contraction 
$F_i\colon \X \ra \X_i$ and a map~$\rho_i\colon Q\rightarrow Q$ fitting into a diagram
\[\begin{array}{rcccl}
&& \X\times_{\bar \X_i, \rho_i}\X && \\
& \swarrow & & \searrow \\
\X &&\stackrel{\theta_{\rho_i}}\dashrightarrow&& \X \\
& \searrow & & \swarrow \\
\downarrow &&\bar \X_i && \downarrow \\
& \swarrow & & \searrow \\
Q &&\stackrel{\rho_i}\longrightarrow && Q 
\end{array}\]
which is just the diagram from Proposition~\ref{many_defs_of_X}(iii)
re-drawn and completed to a fiber product on the top.  
Let \[U_i = \O_{\X\times_{\bar \X_i, \rho_i}\X}\in\db{\X\times_{\rho_i}\X}\]
be the structure sheaf of this fiber product. 
 
\begin{theore} \label{generators} There is a~$W_\Xi$-invariant open subset
$0\in S\subset Q$, such that the restriction of each pair~$(U_i, \rho_i)$ to 
the pullback 
family~$e_S\colon\X_S\ra S$ is in the group of relative equivalences~$\daut{\X_S/S}$.
\end{theore} 

The proof of this theorem relies on the following fact, which will be very 
important also later.

\begin{propositio}\label{flat}
The morphism~$e_i\colon \X\times_{\bar \X_i, \rho_i}\X\rightarrow Q$ is flat and
has reduced fibers. 
\end{propositio} 
\proof Using the notation and constructions of Section~\ref{sec!defs_of_X},
$\X\ra Q$ factors through a morphism~$\X\ra B\times Q$. Hence~$e_i$ factors through 
$\X\times_{\bar \X_i, \rho_i}\X\rightarrow B\times Q$ which is the 
pullback of a morphism 
$\tilde\X\times_{\tilde\X_i, \rho_i}\tilde\X\rightarrow \sH$ along the
natural map~$B\times Q\ra \sH$. It is clearly enough to prove that this latter 
morphism is flat with reduced fibers. However, this statement is (\'etale) local 
over the base, hence it suffices to prove it over the \'etale open set 
$\h_\Delta\times B_l$ of~$\sH$; recall that  
$\{B_l\}$ is an \'etale open covering of the curve~$B$. Over
$\h_\Delta\times B_l$, everything is a pullback along the map 
$\h_\Delta\times B_l\ra \h_\Delta$, so finally it is enough 
to show that the morphism  
\[d_I\colon \Y\times_{\bar\Y_I, r_I}\Y\rightarrow \h_\Delta\]
has the stated properties. Here~$I$ is the~$A$-orbit of nodes of~$\Delta$
corresponding to the node~$i$ of~$\Xi$, and~$r_I$ is the 
corresponding~$A$-fixed element of~$W_\Delta$. 

Since the morphism~$d_I$ is surjective with smooth target~$\h_\Delta$, 
by~\cite[15.2.3 and Remark (v)]{ega}, using also~\cite[14.4.2]{ega}, 
it is flat once its fibers are reduced and its domain 
$\Y\times_{\bar\Y_I, r_I}\Y$
irreducible and equidimensional over~$\h_\Delta$. Equidimensionality is clear, 
so the issue is to prove that the fibers are reduced and the domain 
irreducible. 

Assume first that~$I=i$ is a single node of~$\Delta$. Then the 
central fiber~$\bar Y_I$ of~$\bar\Y_I\ra \h_\Delta$ is a surface with an 
ordinary surface double point, and 
the total family is a deformation family of this surface, where the double point
survives on a codimension one subspace~$\Pi\subset \h_\Delta$. The map
$r_i\colon \h_\Delta\ra \h_\Delta$ is the reflection in~$\Pi$. 
Finally the family~$\Y\rightarrow \h_\Delta$
is a simultaneous resolution of~$\bar\Y_I\ra \h_\Delta$, constructed simply 
by blowing up the singular locus. Hence near the singularity, up to a local 
analytic change in coordinates I~can simply write 
\[\begin{array}{ccccl} 
\bar\Y_i & \cong &  \{xy=z^2-t_1^2\} & \subset& \A^{n+3}_{x,y,z,t_1,\ldots, t_n}  \\
\downarrow&&&&\downarrow \\
\h_\Delta & &=&  & \A^n_{t_1,\ldots, t_n}
\end{array}
\]
with singular locus~${\rm Sing}(\bar\Y_i)=\{x=y=z=t_1=0\}$ mapping to 
$\Pi=\{t_1=0\}\subset \h_\Delta$, the fixed locus of~$r_i\colon t_1\mapsto -t_1$.
The resolution~$\Y$ can be constructed explicitly as the graph of the rational 
map~$\bar\Y_i\dashrightarrow\P^1$ defined by  
$(x,y,z,t_1,\ldots, t_n)\mapsto (x:(z-t_1))$. 
Using the affine variable~$s=x/(z-t_1)$, one affine piece of this graph is
\[\Y^{(1)}\cong\{ys= z+t_1\}\subset\A^{n+3}_{y,z,s,t_1,\ldots, t_n}.\]
Hence the fiber product has an affine open piece
\[(\Y\times_{\bar\Y_i, r_i}\Y)^{(1,1)}= \left\{\begin{array}{rcl} ys_1 & = & z+t_1\\ys_2 & = & z-t_1 \\ s_2(z+t_1) & = & s_1(z-t_1) \end{array}\right\}\subset\A^{n+4}_{y,z,s_1, s_2,t_1,\ldots, t_n}\]
which is isomorphic to the hypersurface
\begin{equation} \Big\{ y(s_1-s_2) = 2t_1\Big\}\subset\A^{n+3}_{y,s_1, s_2,t_1,\ldots, t_n}.\label{irred}\end{equation}
The map to~$\h_\Delta$ is still given by projection to the~$t_i$ coordinates. 
The equation in~(\ref{irred}), together with similar equations for the other 
affine pieces, show that~$\Y\times_{\bar\Y_i, r_i}\Y$ is irreducible, and the map 
to~$\h_\Delta$ has reduced fibers. This concludes the proof for the case 
when~$i=I$ is a single node of~$\Delta$. The other cases reduce to this, since 
locally the morphism $Y\rightarrow \bar Y_I$ contracts a union of disjoint 
rational curves to ordinary double points. 
\ethrm

\vspace{.05in}\noindent {\sc Proof of Theorem~\ref{generators}} \hspace{.05in}
By Proposition~\ref{fiberwise}, it is enough to show
that the fiberwise restricted kernels  
\[U_{i,s} = \pull y_{s}^*U_i\in\db{X_s\times X_{\rho_i(s)}}\]
are invertible, 
where~$y_s\colon X_s\times X_{\rho_i(s)}\hookrightarrow \X\times_{\rho_i}\X$
is fiber inclusion. 
By Proposition~\ref{flat}, 
$U_i$ is a flat family of structure sheaves over~$Q$, and 
hence the derived restriction~$\pull y_{s}^*U_i$ is isomorphic to the ordinary 
restriction~$y_{s}^*U_i$, which in turn is isomorphic to the structure sheaf
$\O_{X_s\times_{\bar X_s}X_{\rho_i(s)}}$. The statement that this sheaf defines
an invertible kernel is already contained in the literature. 
There are two cases. 
 
Suppose first that~$s\in Q\cap T_i$. Then by Proposition~\ref{defs_of_X}(iii), 
$s$ is a fixed point of~$\rho_i$, hence 
\[U_{i,s}\cong\O_{X_s\times_{\bar X_{i,s}}X_s} \in\db{X_s\times X_s}.\] 
On the other hand, the 
contraction~$f_{i,s}\colon X_s\rightarrow \bar X_{i,s}$ contracts a single ruled 
surface~$D_{i,s}$ inside~$X_s$ to a smooth curve~$B_{i,s}$ 
(which is either~$B$ or~$\tilde B$). There is an exact sequence of sheaves
on~$X_s\times X_s$ 
\[ 0\ra U_{i,s} \ra \O_{\Delta_{X_s}}\oplus \O_{D_{i,s}\times_{B_{i,s}}D_{i,s} }\ra \O_{\Delta_{D_{i,s}}}\ra 0, 
\]
where~$\Delta_{X_s}$ and~$\Delta_{D_{i,s}}$ are respective diagonals in 
$X_s\times X_s$. Hence the kernel~$U_{i,s}$ is isomorphic to the kernel 
\begin{equation} \label{ker_h}\left( \O_{\Delta_{X_s}}\oplus \O_{D_{i,s}\times_{B_{i,s}}D_{i,s} }\ra \O_{\Delta_{D_{i,s}}}\right)\in \db{X_s\times X_s}. 
\end{equation} 
This kernel was introduced in~\cite[(4.31)]{horja1} and its invertibility
proved in~\cite[Theorem 2.9 and Remark 2.12]{horja2}. 

Next suppose that~$s\in Q \cap (T\setminus T_i)$. Then by 
Proposition~\ref{many_defs_of_X}(iii), the birational map 
\[\theta_{\rho_i,s}\colon X_s\dashrightarrow  X_{\rho_i(s)}\] is a flop of a 
disjoint union of smooth rational curves with normal bundle 
$\O_{\P^1}(-1,-1)$ or $\O_{\P^1}(0,-2)$. The kernel~$U_s$ is the structure sheaf 
of the graph of this flop. This kernel was shown to be invertible
in~\cite[Theorem 3.6 and Remark]{bo}. 
\ethrm

\begin{remar}\rm I~sketch an alternative proof of the invertibility of
the kernel~$U_{i,s}$, which avoids a case division and also throws some light 
on the origin of this kernel. The claim is that~$U_{i,s}$ is the 
universal perverse coherent point sheaf on~$X_s\times X_{r_i(s)}$ with respect 
to the contraction~$f_{i,s}\colon X_s\ra\bar X_{i,s}$, and hence it is invertible.
In particular, the variety~$X_{r_i(s)}$ is the fine moduli space of perverse
point sheaves on~$X_s$ for the contraction $f_{i,s}$. 
Here I~am using the terminology
of~\cite{bridgeland_flops}; the essential point is that~$f_{i,s}$ has fibers
of dimension at most one, so Bridgeland's theory applies. The proof of the claim is
not very difficult given the machinery of~\cite{bridgeland_flops}.
\end{remar} 

For purposes of brevity I~will denote the family over the open set~$S\subset Q$ 
also by~$\X\ra S$, and restrict all contractions, 
the~$W_\Xi$-action and the relative kernels~$(U_i, \rho_i)$ 
to this family without further notice. The properties spelled out 
in Propositions~\ref{defs_of_X}--\ref{many_defs_of_X} continue to hold;
the latter of course under the assumption that~$S$ contains 
a sufficiently general point of~$T$. 

\subsection{The main results}\label{sub!main} The first main result of the paper
is that the derived category of the threefold~$X$ and 
that of its deformation space carry an action of an Artin group. 

\begin{theore} Let~$\X\rightarrow S$ be a finite-dimensional deformation 
space of the threefold~$X$ satisfying the conclusion of Theorem~\ref{generators}. 
Then there are homomorphisms
\begin{equation}\label{main_hom1} 
B_\Xi \longrightarrow \daut{\X/S} 
\end{equation}
and
\begin{equation}\label{main_hom2}
B_\Xi \longrightarrow \daut{X}. 
\end{equation}
\label{main_theorem}\end{theore}
\proof Define the map~(\ref{main_hom1}) by mapping the generators 
$\rho_i$ of the Artin group~$B_\Xi$ to the element 
$(U_i, \rho_i)$ of Theorem~\ref{generators}. Since~$W_\Xi$ 
fixes~$0\in S$, 
I~can define~(\ref{main_hom2}) by restricting these kernels to the central 
fiber. The point is to prove that the braid relations of~(\ref{braidrelations}) 
defining~$B_\Xi$ are satisfied for these kernels. Since (derived) restriction
commutes with kernel composition in smooth families, it is enough to show
that~(\ref{main_hom1}) is a group homomorphism.

Take a pair of nodes~$(i,j)$ of the Dynkin diagram~$\Xi$. Set
\[
(V_l, \phi_l)=\underbrace{(U_i, \rho_i)\circ (U_j, \rho_j)\circ\ldots}_{m_{ij}} 
\]
for the left hand side of the braid relation of~(\ref{braidrelations}) for the
pair of nodes~$i, j$ and similarly~$(V_r, \phi_r)$ for the right hand side.
One part is easy: the automorphisms~$\rho_i$ and~$\rho_j$ of the base~$S$ satisfy 
the relations of the Coxeter group~$W_\Xi$, and consequently also the braid 
relation; hence~$\phi_l=\phi_r$ which I~will denote simply by~$\phi$. 

Next note that for~$k=i,j$, the sheaf~$U_k$ is the structure sheaf of a subscheme 
of~$\X\times_{\rho_k}\X$, in other words there is a surjective morphism 
$\O_{\X\times_{\rho_k}\X}\rightarrow U_k$. 
By a repeated use of Lemma~\ref{arrows_from_O}, this implies that there are
induced arrows
\begin{equation} \label{arrow} 
\begin{array}{rcccl} && \O_{\X\times_\phi\X}\\
& \swarrow && \searrow \\
V_l&& && V_r
\end{array}
\end{equation} 
in~$\db{\X\times_{\phi}\X}$.

To continue, assume that~$m_{ij}=3$; this will only simplify notation, the other 
cases being identical. Let~$\rho_{ij}=\rho_j\circ\rho_i$. Let also 
\[\tilde\X= \X\times_{\rho_i}(\X\times_{\rho_j}(\X\times_{\rho_i}\X))\]
with maps~$p_{12}\colon \tilde\X\ra\X\times_{\rho_i}\X$ etc. For~$s\in S$, let 
\[x_s\colon\tilde X_s= X_s\times X_{\rho_i(s)}\times X_{\rho_{ij}(s)}\times X_{\phi(s)}\hookrightarrow \tilde\X\] and 
\[y_s\colon X_s\times X_{\phi(s)}\hookrightarrow\O_{\X\times_\phi\X}\]
be inclusion maps of fibers, with projection maps 
$p_{14s}\colon\tilde X_s\ra X_s\times X_{\phi(s)}$ etc. 
The derived restriction of diagram~(\ref{arrow}) is a diagram 
\begin{equation} \label{arrow_s}
\begin{array}{rcccl} && \O_{X_s\times X_{\phi(s)}}\\
& \swarrow && \searrow \\
V_{l,s}&& && V_{r,s}
\end{array}
\end{equation} 
in~$\db{X_s\times X_{\phi(s)}}$. Now compute: 
 
\begin{equation}\label{vls}\begin{array}{rcl} V_{l,s} & \cong & \pull y_s^*\left(\push p_{14*}(p_{34}^*(U_i)\tensor p_{23}^*(U_j)\tensor p_{12}^*(U_i))\right) \\
&\cong& \push p_{14s*}\left(\pull x_s^*(p_{34}^*(U_i)\tensor p_{23}^*(U_j)\tensor p_{12}^*(U_i))\right) \\
&\cong& \push p_{14s*}\left(
\O_{X_s\times X_{\rho_i(s)}\times X_{\rho_{ij}(s)}\times_{\bar X_{i,\rho_{ij}(s)}} X_{\phi(s)}}\right.
\\
&&\tensor
\O_{X_s\times X_{\rho_i(s)}\times_{\bar X_{j,\rho_{i}(s)}} X_{\rho_{ij}(s)}\times X_{\phi(s)}}
\\
&&\left.\tensor
\O_{X_s\times_{\bar X_{j,s}} X_{\rho_i(s)}\times X_{\rho_{ij}(s)}\times X_{\phi(s)}}
\right) . 
\end{array}\end{equation}
The first isomorphism uses Lemma~\ref{three_ker}, the second 
follows from a slight generalization of~\cite[Lemma 1.3]{bo} to the 
quasiprojective case, and the last uses 
the flatness result Proposition~\ref{flat}. There is a similar 
computation for~$V_{r,s}$. 

I~now distinguish two cases. First assume that~$s\in S$ is sufficiently
general. It is easy to see that the three subschemes of~$\tilde X_s$ 
appearing in the last expression of~(\ref{vls})
are transversal, so the (derived) tensor product of their 
intersections is isomorphic in $\db{\tilde X_s}$ to
the structure sheaf of their intersection
\[\tilde C_{l,s}= X_s\times_{\bar X_{i,s}}  X_{\rho_i(s)}\times_{\bar X_{j,\rho_{i}(s)}} X_{\rho_{ij}(s)}\times_{\bar X_{i,\rho_{ij}(s)}} X_{\phi(s)}\]
in~$\tilde X_s$. To understand this intersection, consider the diagram
\begin{equation}\label{cool}\begin{array}{ccccccccccccc}
X_s  &&\dashrightarrow &&  X_{\rho_i(s)}  &&\dashrightarrow&& X_{\rho_{ij}(s)} &&\dashrightarrow&&  X_{\phi(s)}\\
&\searrow && \swarrow &&\searrow && \swarrow&&\searrow && \swarrow\\
&&\bar X_{i,s} &&&& \bar X_{j,\rho_{i}(s)} &&&& \bar X_{i,\rho_{ij}(s)}
\end{array}\end{equation}
The support of~$\tilde C_{l,s}$ in~$\tilde X_s$ is the set of quadruples
$(p_1,\ldots, p_4)\in\tilde X_s$ such that~$p_i, p_{i+1}$ have the same images
under appropriate arrows in the diagram~(\ref{cool}). 

Since~$s$ is sufficiently general, the variety~$X_s$ contains a disjoint union 
of~$(-1,-1)$ 
curves indexed by positive roots. The chain of birational maps in~(\ref{cool}) 
flops, consequtively, the disjoint rational curves on~$X_s$ indexed by the 
positive roots~$\mu_i, \mu_i+\mu_j, \mu_j\in\Sigma_\Xi^+$. In other words, the 
composition of the three flops, the map~$X_s\dashrightarrow X_{\phi(s)}$, is 
the flop of the disjoint set of all these curves. Hence a quadruple 
$(p_1,\ldots, p_4)$ is determined by the pair~$(p_1, p_4)$, and in particular
$p_{14s}$ restricted to the reduced subscheme~$\tilde C_{l,s}$ of~$\tilde X_s$ is 
an isomorphism onto its image in~$X_s\times X_{\phi(s)}$. 

On the other hand, the set of positive roots~$\{\mu_i, \mu_i+\mu_j, \mu_j\}$ 
is exactly the set of all positive roots of the sub-root 
system of type~$A_2$ of the root system of~$\Xi$ spanned by 
the nodes~$(i,j)$ (remember~$m_{ij}=3$). The reflection~$\phi\in W_\Xi$ maps 
exactly these roots to negative roots. So the map 
$X_s\dashrightarrow X_{\phi(s)}$ factors as 
$X_s\ra\bar X_{ij,s} \leftarrow X_{\phi(s)}$ in the notation of Proposition~\ref{many_defs_of_X}(iiib). 
Hence the image of~$\tilde C_{l,s}$ under~$p_{14s}$ is the reduced subscheme 
\[C_{l,s}=X_s\times_{\bar X_{ij,s}} X_{\phi(s)}\hookrightarrow X_s\times X_{\phi(s)}.\]
So finally
\[ V_{l,s}\cong p_{14s*}\O_{\tilde C_{l,s}}\cong \O_{C_{l,s}}
\]
and the argument also shows that the map on the left hand side of 
diagram~(\ref{arrow_s}) is just the natural surjection 
$\O_{X_s\times X_{\phi(s)}}\ra\O_{C_{l,s}}$. 

Repeating this argument also for~$V_{r,s}$, I~obtain that in this case the 
diagram~(\ref{arrow_s}) is isomorphic to the diagram 
\begin{equation} \label{arrow_gen}
\begin{array}{rcccl} && \O_{X_s\times X_{\phi(s)}} \\
& \swarrow && \searrow \\
\O_{X_s\times_{\bar X_{ij,s}} X_{\phi(s)}}&& \stackrel{\sim}{\longrightarrow}&& \O_{X_s\times_{\bar X_{ij,s}} X_{\phi(s)}}
\end{array}
\end{equation} 
of sheaves on~$X_s\times X_{\phi(s)}$, with the vertical arrows being surjective.

The other case is when~$s\in S$ is not sufficiently general. Transversality
of the subschemes in~(\ref{vls}) still holds, hence the tensor product in the last 
expression of~(\ref{vls}) is isomophic to the structure sheaf of a 
subscheme~$\tilde C_{l,s}$ of~$\tilde X_s$, a correspondence subscheme with respect 
to the diagram~(\ref{cool}). It is no longer true that the horizontal birational
maps in~(\ref{cool}) modify~$X_s$ on a set of disjoint loci, but in any case
$p_{14s}$ restricted to~$\tilde C_{l,s}$ factors as 
\[p_{14s}|_{\tilde C_{l,s}}\colon \tilde C_{l,s}\twoheadrightarrow C_{l,s}\hookrightarrow X_s\times X_{\phi(s)}.\]
The fibers of~$p_{14s*}|_{\tilde C_{l,s}}$ are quadruples~$(p_1,\ldots, p_4)$ 
mapping to a given pair~$(p_1,p_4)\in X_s\times X_{\phi(s)}$ and~$p_i, p_{i+1}$ 
mapping to the same image under the appropriate maps in~(\ref{cool}).
There are several possible configurations, depending on~$s$ and the exceptional 
loci; for example if~$s\in T_i\cap T_j$ then one possibility is 
that~$p_1$,~$p_2=p_3$ and~$p_4$ all lie in the same fiber of~$f_{i,s}$. However, 
it is easy to check that in all cases when the reduced fiber of 
$p_{14s*}|_{\tilde C_{l,s}}$ is not a point, 
it is isomorphic to a fiber of either~$f_{i,s}$ or 
$f_{j,s}$, in other words to a projective line. Hence by Lemma~\ref{pushdown}, 
\[ V_{l,s}\cong\push p_{14s*}\O_{\tilde C_{l,s}} \cong \O_{C_{l,s}}\in\db{X_s\times X_{\phi(s)}}.\] 
Using the same reasoning also for~$V_{r,s}$, the diagram~(\ref{arrow_s}) is 
isomorphic to a diagram of sheaves on~$X_s\times X_{\phi(s)}$
\begin{equation}\label{arrow_spec}
\begin{array}{rcccl} && \O_{X_s\times X_{\phi(s)}}\\
& \swarrow && \searrow \\
\O_{C_{l,s}}&& && \O_{C_{r,s}}
\end{array}
\end{equation} 
with surjective arrows. 
 
Diagrams~(\ref{arrow_gen}) and~(\ref{arrow_spec}) imply 
by~\cite[Lemma 4.3]{bridgeland} that~$V_l$ and~$V_r$ are 
sheaves on~$\X\times_\phi\X$, flat over~$S$. Moreover, since pullback is right 
exact, the arrows in diagram~(\ref{arrow}) are necessarily 
surjective maps of sheaves; in other words,~$V_l$ and~$V_r$ are structure 
sheaves of subschemes of~$\X\times_\phi\X$. If I~further assume that 

\vspace{0.1in}

\noindent{$(\star)$}~$\M$ is a moving linear system on~$B$, and 
the finite-dimensional family~$e\colon \X\ra S$ contains sufficiently general 
deformations of~$X$,

\vspace{0.1in}

\noindent then sufficiently general points form an open dense subset of 
$S$. Hence Lemma~\ref{isom_for_subshemes}, together with~(\ref{arrow_gen}), 
allows me to conclude that~(\ref{arrow}) can be extended to a diagram 
\[
\begin{array}{rcccl} && \O_{\X\times_\phi\X}\\
& \swarrow && \searrow \\
V_l&&\stackrel{\sim}{\longrightarrow} && V_r. 
\end{array}
\]
So composition of the relative kernels~$(U_i, \rho_i)$ and
$(U_j, \rho_j)$ in the two different ways gives isomorphic relative kernels; 
hence, assuming {$(\star)$}, the braid relation holds up to isomorphism. 

To remove assumption {$(\star)$}, let~$B=\cup B_{\beta}$ be a
finite decomposition into quasiprojective (e.g.~affine) curves so that 
the restriction~$\M_\beta= \M|_{B_\beta}$ moves on~$B_\beta$. For every~$\beta$, 
let~$\pi_\beta\colon X_\beta\ra B_\beta$ be the restriction of~$\pi$ over~$B_\beta$, and 
let~$e_\beta\colon\X_\beta\ra H^0(B_\beta, \M_\beta)$ be
the family of deformations of~$X_\beta$ constructed in Proposition~\ref{defs_of_X}. 
There is a~$W_\Xi$-equivariant natural injection
\[ n_\beta\colon H^0(B, \M)\hookrightarrow H^0(B_\beta, \M_\beta). 
\] 
Let~$S_\beta$ be a~$W_\Xi$-invariant finite-dimensional 
subspace of~$H^0(B_\beta, \M_\beta)$ containing both the image of~$S$ under 
the injection~$n_\beta$ and a sufficiently general point of 
$H^0(B_\beta, \M_\beta)$. There is an induced family  
$e_\beta\colon\X_\beta\ra S_\beta$ with central fiber~$X_\beta$. 

By construction, the family~$e_\beta\colon\X_\beta\ra S_\beta$ satisfies 
assumption~$(\star)$ for each~$\beta$. On the other hand, the natural injection 
$n_\beta|_S\colon S\hookrightarrow S_\beta$ is~$W_\Xi$-equivariant by construction, 
so if I~restrict to families 
$e_{\beta, S}\colon \X_{\beta,S}\rightarrow S$ then the above discussion applies to 
these families. In particular, the restriction of 
diagram~(\ref{arrow}) to~$\X_{\beta, S}\times_\phi\X_{\beta,S}$ can be extended 
to a diagram of sheaves 
\begin{equation} \label{beta_arrow} 
\begin{array}{rcccl} && \O_{\X_{\beta,S}\times_\phi\X_{\beta,S}}\\
& \swarrow && \searrow \\
V_{l,\beta, S}&& \stackrel{\sim}{\longrightarrow} && V_{r,\beta,S} 
\end{array}\end{equation}
with the vertical maps being surjective. The horizontal isomorphisms
in~(\ref{beta_arrow}) are compatible with surjections from a fixed sheaf, 
so they can be glued to an isomorphism
\begin{equation} 
\begin{array}{rcccl} && \O_{\X\times_\phi\X}\\
& \swarrow && \searrow \\
V_l&&\stackrel{\sim}{\longrightarrow} && V_r
\end{array}
\end{equation} 
of sheaves on~$\X\times_\phi\X$ extending~(\ref{arrow}). Hence the braid 
relation holds between~$(U_i, \rho_i)$ and~$(U_j, \rho_j)$ 
with no extra assumption. 

The proofs in the cases~$m_{ij}=4,6$ are, up to writing out longer
expressions, identical. The proof for~$m_{ij}=2$ is in fact easier, 
since the exceptional loci of~$f_{i,s}$ and~$f_{j,s}$ are disjoint for 
all~$s$, hence there is no need for a case distinction and assumption {$(\star)$}. 
These cases correspond to sub-digrams of~$\Xi$ type~$B_2$, $G_2$ 
and~$A_1\times A_1$ respectively. The proof of Theorem~\ref{main_theorem} 
is complete. 
\ethrm

The next result shows that in certain cases, the Artin group action 
on the derived category of~$X$ is faithful. 

\begin{theore} Assume that the diagram~$\Xi$ describing the configuration of 
exceptional surfaces in~$X$ is of type~$A_n$ or~$C_n$. Then 
the map~$B_\Xi\longrightarrow \daut{X}$ is injective. 
\label{faith}\end{theore} 
\proof Take any point~$b\in B$ and let~$j_b\colon Y_b\rightarrow X$ be the fiber of 
$\pi\colon X\ra B$ over~$b$. Rational exceptional 
curves~$E_{j,b}$ in the surface~$Y_b$ are indexed by nodes~$j$ of the 
simply laced Dynkin diagram~$\Delta$ lying over~$\Xi$. 

Let~$k_b\colon Y_b\times Y_b\ra X\times X$; for nodes~$i$ of~$\Xi$, let 
$U_{i,b}=\pull k_b^*(U_{i,0})$ be the restriction to~$Y_b\times Y_b$ 
of the kernel~$U_{i,0}$ (note~$X=X_0$ for~$0\in S$). 
By standard arguments, there is a commutative diagram of functors
\[ \begin{array}{ccc} \db{X} & \stackrel{\pull k_b^*}{\longrightarrow} & \db{Y_b}\\
\mapdownright{\Psi^{U_{i,0}}} & & \mapdownright{\Psi^{U_{i,b}}}\\
\db{X} & \stackrel{\pull k_b^*}{\longrightarrow} & \db{Y_b}
\end{array}\] 
and maps \[B_\Xi\rightarrow \daut{X}\rightarrow \daut{Y_b}\]
where the second arrow is restriction to the fiber over~$b$.

Suppose now that~$\Xi$ is of type~$A_n$. 
Then the following holds (for a proof, see below):
\begin{lemm} The kernel~$U_{i,b}\in\db{Y_b\times Y_b}$ is 
isomorphic to the kernel defining the (inverse) twist functor 
$T'_{\E_i}$ of~\cite{st} for the sheaf~$\E_i=\O_{E_{i,b}}(-1)$, where 
$E_{i,b}\subset Y_b$ is the exceptional rational curve corresponding to the node 
$i$ of~$\Delta$.   
\label{st1}\end{lemm} 
The map~$B_\Xi\rightarrow \daut{Y_b}$ defined by mapping the Artin group 
generators to the twist functors~$T_{\E_i}$ is injective 
by~\cite[Theorem 2.18]{st}. Hence the map~$B_\Xi\rightarrow \daut{X}$ 
must be injective as well. 

If~$\Xi$ is of type~$C_n$, then it has two kinds of nodes: 
one representing a single node of the simply laced diagram~$\Delta$, and the 
others representing an orbit~$\{i_1, i_2\}$ of nodes. 
For the first type of node, Lemma~\ref{st1} continues to hold; for the second, 
it gets replaced by 
\begin{lemm} The kernel~$U_{i,b}\in\db{Y_b\times Y_b}$ is 
the composite of the commuting kernels defining the (inverse) 
twist functors~$T'_{\E_{i_1}}$ and~$T'_{\E_{i_2}}$. \label{st2}\end{lemm} 
\noindent Hence, recalling the proof of Lemma~\ref{maps_braidgp},
in this case there is a commutative diagram 
\[ \begin{array}{ccc} B_\Xi & \longrightarrow & \daut{X}\\
\Big\downarrow & & \Big\downarrow\\
B_\Delta & \longrightarrow & \daut{Y_b}. 
\end{array}\] 
The bottom horizontal arrow is injective by~\cite[Theorem 2.18]{st} again; 
the left hand vertical arrow is injective by Lemma~\ref{maps_braidgp}.
Hence the composite is injective; so~$B_\Xi\rightarrow \daut{X}$
must be injective as well. 
\ethrm

\vspace{.05in}\noindent {\sc Proof of Lemmas~\ref{st1} and \ref{st2}} \hspace{.05in}
Suppose first that~$A$ is trivial. Since all sheaves appearing 
in~(\ref{ker_h}) are flat with respect to the projection~$X\times X\ra B$, the 
kernel~$U_{i,b}$ on~$Y$ is isomorphic to the kernel 
\[\left( \O_{\Delta_{Y_b}}\oplus \O_{E_{i,b}\times E_{i,b} }\ra \O_{\Delta_{E_{i,b}}}\right)\in \db{Y_b\times Y_b} \]
where recall~$E_{i,b}\subset Y_b$ is one of the exceptional rational curves. 
Let~$y\colon E_{i,b}\hookrightarrow Y_b$ denote the inclusion of the rational 
curve in the surface, and let 
$x\colon E_{i,b}\times E_{i,b}\hookrightarrow Y_b\times Y_b$ 
be the induced inclusion. Then up to isomorphism in~$\db{Y_b\times Y_b}$,
\begin{eqnarray*}
U_{i,b}&\cong & \cone\Big\{\O_{\Delta_{Y_b}}\ra \Big(x_*\O_{E_{i,b}\times E_{i,b}}\ra x_*\O_{\Delta_{E_{i,b}}}\Big)\Big\}\\
& \cong & \cone\Big\{\O_{\Delta_{Y_b}}\ra x_*\Big(\O_{E_{i,b}\times E_{i,b}}(-1,-1)\Big)\Big\}\\
& \cong & \cone\left\{\O_{\Delta_{Y_b}}\ra \right. \\
& & \left.x_*\left(q_2^*\Big(\rhom_{\O_{E_{i,b}}}(\O_{E_{i,b}}(-1),\omega_{E_{i,b}})\Big)\tensor q_1^*\Big(\O_{E_{i,b}}(-1)\Big)\right)\right\}
\end{eqnarray*}
where~$q_1, q_2\colon E_{i,b}\times E_{i,b}\ra E_{i,b}$ are projections. 
It is easiest to conclude now
using results of~\cite{horja2}. Comparing the above expression 
with~\cite[(2.7) and (2.23)]{horja2} shows that~$U_{i,b}$ is isomorphic
to the invertible kernel on~$Y_b$ defined by the relatively spherical 
sheaf~$\O_{E_{i,b}}(-1)$ with respect to the diagram 
\[\begin{array}{ccc}E_{i,b} & \stackrel{y}\hookrightarrow & Y_b \\ \Big\downarrow \\ \star\end{array}\] 
a kernel which by~\cite[Example 4.1]{horja2} is isomorphic to the kernel which 
gives rise to the (inverse) twist functor of~\cite{st} defined by the spherical 
sheaf~$y_*\O_{E_{i,b}}(-1)$. This proves 
Lemma~\ref{st1}. Lemma~\ref{st2} follows also on noting that in that case the 
contraction~$f_i$ restricts to~$Y_b$ as the contraction of two disjoint 
exceptional curves. 
\ethrm

\subsection{Projective examples}\label{sub!proj}

Let~$\bar X$ be a projective threefold with a curve 
of singularities
\[B={\rm Sing}(X)\hookrightarrow \bar X,\] 
such that along the
curve,~$\bar X$ has compound du Val singularities of uniform~$ADE$ type. 
The iterated blowup of the singular locus~$f\colon X\ra \bar X$ is a resolution of 
singularities, cf.~\cite{pagoda}, and the exceptional locus consists of a set of 
geometrically ruled surfaces~$\{\pi_j\colon D_j\ra B_j\}$ intersecting in one of the 
configurations~$\Xi=\Delta/A$ described in Section~\ref{threefold}.

\begin{theore} Assume as usual that~$(\Delta,A)\neq (A_{2n}, \Z/2)$.
The derived category~$\db{X}$ carries an action of the 
Artin group~$B_\Xi$. In case~$\Xi$ is of type~$A_n$ or~$C_n$, this action is 
faithful.
\end{theore} 
\proof For~$j$ a node of the diagram~$\Xi$, define a kernel~$U_j$ on~$X$ by
\[U_j=\left( \O_{\Delta_X}\oplus \O_{D_j\times_{B_j}D_j}\ra \O_{\Delta_{D_j}}\right)\in \db{X\times X}; 
\]
this is just the kernel of~\cite[(4.31)]{horja1}, proved to be invertible 
in~\cite{horja2}. The point is that this definition makes sense whether or not 
there is a contraction morphism on~$X$ contracting~$D_j$ alone. Define the map 
\[ B_\Xi\longrightarrow \daut{X}\]
by mapping the generator~$R_j$ of~$B_\Xi$ to the kernel~$U_j$ on~$X$. The issue 
is again to prove the braid relations. As before, take a pair of nodes~$(i,j)$
of~$\Xi$ and let~$V_l, V_r\in\db{X\times X}$ be the composite kernels on the 
two sides of the braid relation for the pair~$(i,j)$.
Note that the interesting part of 
the computation of all these kernels takes place in an \'etale neighbourhood of 
the exceptional set; away from such a neighbourhood,~$V_l$ and~$V_r$ are 
obviously isomorphic to the structure sheaf of the diagonal. There is an
\'etale open covering of a neighbourhood of~$B\subset\bar X$, such that 
on the inverse image of this covering on~$X$ the restrictions of~$V_l$ 
and~$V_r$ are isomorphic by the proof of Theorem~\ref{main_theorem}, these
isomorphisms being compatible on intersections. Hence by descent, there is an 
isomorphism~$V_l\cong V_r$ on~$X\times X$, and so the braid relations hold
up to isomorphism. 

To prove faithfulness, argue as in the proof of 
Theorem~\ref{faith}: take a quasiprojective surface~$\bar Y_s\subset\bar X$ 
intersecting the singular locus~$B\subset\bar X$ transversally at 
$p\in B$ and in no other points, let~$Y_s\subset X$ be its resolution, 
restrict the kernels~$U_i$ to~$Y_s$ using Lemmas~\ref{st1} and \ref{st2} and 
appeal to the faithfulness result of~\cite{st}. 
\ethrm

Examples of varieties~$\bar X$ with a curve of singularities of uniform 
type~$A_n$ can be found among hypersurfaces or complete intersections in weighted 
projectice spaces; compare for example~\cite{kmp}. The resolution~$X$ is then 
embedded in a (partial) resolution of the ambient space, typically with 
$n$ distinct divisors over the relevant singular locus; hence the configuration
in~$X$ is still of type~$A_n$. Such varieties can be found and in low
codimension classified using the graded ring method pioneered by Reid; 
see the (from the present point of view not very interesting)~$A_1$ case 
in~\cite{sz_wci} and the general case in~\cite{anita}.  
Examples of type~$(A_n, \Z/2)$ can be constructed as quotients; 
see~\cite[Examples 4.3]{eg} for an explicit example. In favourable cases, 
the local deformations described in Proposition~\ref{many_defs_of_X} are realized 
as actual projective deformations. In such cases, the action of the 
Artin group on~$\db{X}$ can be extended to an action 
by relative equivalences over its local universal family, in an analogous way 
to the statement of Theorem~\ref{main_theorem}.
I~leave it to the reader to formulate the precise statement. 

\begin{remar} \rm The action of the Artin group $B_\Xi$ on the derived category 
of the threefold~$X$ gives rise 
to actions on even and odd cohomology, using the Chern class map. In the case
when~$X$ has trivial canonical bundle,~$H^{2,1}(X)\cong H^1(X,\Theta_X)$ is a 
direct summand of odd cohomology, and it is preserved by the action. 
Hence the braid group acts on the tangent space to the deformation space, and 
it is easy to see that this action factors through the reflection group~$W_\Xi$. 
The action on even cohomology can in turn be restricted to the Picard group 
to get an action of~$W_\Xi$ there. 
Some of these actions were known before; e.g.~\cite{wi_ell} discusses 
the case of elliptic ruled surfaces, whereas~\cite{kmp} has a 
symmetric group action in the case of Type~$A$. The action of the Artin group on 
the derived category shows the uniform origin of all these actions. 
\end{remar}

\vspace{0.2in} 

\noindent {\small \sc Department of Mathematics, Utrecht University, PO.~Box 80010, NL-3508 TA Utrecht, The Netherlands} 

and 

\noindent {\small \sc Alfr\'ed R\'enyi Institute of Mathematics, Hungarian Academy of Sciences, PO.~Box 127, H-1364 Budapest, Hungary}

\noindent {\small E-mail address: \tt szendroi@math.uu.nl}
\end{document}

%% file: conf.pstex_t
\begin{picture}(0,0)%
\includegraphics{conf.pstex}%
\end{picture}%
\setlength{\unitlength}{4144sp}%
\begingroup\makeatletter\ifx\SetFigFont\undefined%
\gdef\SetFigFont#1#2#3#4#5{%
  \reset@font\fontsize{#1}{#2pt}%
  \fontfamily{#3}\fontseries{#4}\fontshape{#5}%
  \selectfont}%
\fi\endgroup%
\begin{picture}(6889,4689)(856,-5020)
\put(3016,-4066){\makebox(0,0)[lb]{\smash{\SetFigFont{12}{14.4}{\familydefault}{\mddefault}{\updefault}{\color[rgb]{0,0,0}$A_3$}%
}}}
\put(856,-3616){\makebox(0,0)[lb]{\smash{\SetFigFont{12}{14.4}{\familydefault}{\mddefault}{\updefault}{\color[rgb]{0,0,0}$\Delta$}%
}}}
\put(856,-4651){\makebox(0,0)[lb]{\smash{\SetFigFont{12}{14.4}{\familydefault}{\mddefault}{\updefault}{\color[rgb]{0,0,0}$\Xi$}%
}}}
\put(6886,-4021){\makebox(0,0)[lb]{\smash{\SetFigFont{12}{14.4}{\familydefault}{\mddefault}{\updefault}{\color[rgb]{0,0,0}$D_4$}%
}}}
\put(4771,-4066){\makebox(0,0)[lb]{\smash{\SetFigFont{12}{14.4}{\familydefault}{\mddefault}{\updefault}{\color[rgb]{0,0,0}$D_4$}%
}}}
\put(6931,-4966){\makebox(0,0)[lb]{\smash{\SetFigFont{12}{14.4}{\familydefault}{\mddefault}{\updefault}{\color[rgb]{0,0,0}$G_2$}%
}}}
\put(1441,-4966){\makebox(0,0)[lb]{\smash{\SetFigFont{12}{14.4}{\familydefault}{\mddefault}{\updefault}{\color[rgb]{0,0,0}$A_2$}%
}}}
\put(1441,-4066){\makebox(0,0)[lb]{\smash{\SetFigFont{12}{14.4}{\familydefault}{\mddefault}{\updefault}{\color[rgb]{0,0,0}$A_2$}%
}}}
\put(3016,-4966){\makebox(0,0)[lb]{\smash{\SetFigFont{12}{14.4}{\familydefault}{\mddefault}{\updefault}{\color[rgb]{0,0,0}$C_2$}%
}}}
\put(4816,-4966){\makebox(0,0)[lb]{\smash{\SetFigFont{12}{14.4}{\familydefault}{\mddefault}{\updefault}{\color[rgb]{0,0,0}$B_3$}%
}}}
\end{picture}

%% file: loci.pstex_t
\begin{picture}(0,0)%
\includegraphics{loci.pstex}%
\end{picture}%
\setlength{\unitlength}{4144sp}%
\begingroup\makeatletter\ifx\SetFigFont\undefined%
\gdef\SetFigFont#1#2#3#4#5{%
  \reset@font\fontsize{#1}{#2pt}%
  \fontfamily{#3}\fontseries{#4}\fontshape{#5}%
  \selectfont}%
\fi\endgroup%
\begin{picture}(5522,5953)(251,-5465)
\put(3196,-736){\makebox(0,0)[lb]{\smash{\SetFigFont{12}{14.4}{\familydefault}{\mddefault}{\updefault}{\color[rgb]{0,0,0}$X$}%
}}}
\put(4951,-466){\makebox(0,0)[lb]{\smash{\SetFigFont{12}{14.4}{\familydefault}{\mddefault}{\updefault}{\color[rgb]{0,0,0}$X_s$}%
}}}
\put(1081,344){\makebox(0,0)[lb]{\smash{\SetFigFont{12}{14.4}{\familydefault}{\mddefault}{\updefault}{\color[rgb]{0,0,0}$X_u$}%
}}}
\put(271,-1726){\makebox(0,0)[lb]{\smash{\SetFigFont{12}{14.4}{\familydefault}{\mddefault}{\updefault}{\color[rgb]{0,0,0}$X_t$}%
}}}
\put(3106,-5416){\makebox(0,0)[lb]{\smash{\SetFigFont{12}{14.4}{\familydefault}{\mddefault}{\updefault}{\color[rgb]{0,0,0}$T$}%
}}}
\put(1036,-5011){\makebox(0,0)[lb]{\smash{\SetFigFont{12}{14.4}{\familydefault}{\mddefault}{\updefault}{\color[rgb]{0,0,0}$t$}%
}}}
\put(4141,-4561){\makebox(0,0)[lb]{\smash{\SetFigFont{12}{14.4}{\familydefault}{\mddefault}{\updefault}{\color[rgb]{0,0,0}$T_2$}%
}}}
\put(2206,-5011){\makebox(0,0)[lb]{\smash{\SetFigFont{12}{14.4}{\familydefault}{\mddefault}{\updefault}{\color[rgb]{0,0,0}$T_1$}%
}}}
\put(1801,-4561){\makebox(0,0)[lb]{\smash{\SetFigFont{12}{14.4}{\familydefault}{\mddefault}{\updefault}{\color[rgb]{0,0,0}$u$}%
}}}
\put(3151,-4561){\makebox(0,0)[lb]{\smash{\SetFigFont{12}{14.4}{\familydefault}{\mddefault}{\updefault}{\color[rgb]{0,0,0}$0$}%
}}}
\put(4816,-4201){\makebox(0,0)[lb]{\smash{\SetFigFont{12}{14.4}{\familydefault}{\mddefault}{\updefault}{\color[rgb]{0,0,0}$s$}%
}}}
\end{picture}